\documentclass{conm-p-l}

\usepackage{amsmath,amssymb,amsfonts}

\usepackage[mathscr]{eucal}
\usepackage{amsthm}

\usepackage{curves}

\usepackage{srcltx}

\usepackage[all]{xy}

\numberwithin{equation}{section}
\renewcommand\qed{{\hspace*{\fill}Q.E.D.\vskip12pt plus 1pt}}
\newcommand\Pic[1]{\hbox{\rm Pic(}#1\hbox{\rm )}}

\newcommand\length{\hbox{\rm length}}

\newcommand\sC{ \mathcal{C}}

\newcommand\sE{ \mathcal{E}}
\newcommand\sF{ \mathcal{F}}
\newcommand\sG{ \mathcal{G}}

\newcommand\sL{ \mathcal{L}}

\newcommand\sO{ \mathcal{O}}

\newcommand\sQ{ \mathcal{Q}}

\newcommand\sT{ \mathcal{T}}

\newcommand\scc{{\mathscr C}}
\newcommand\scr{{\mathscr R}}
\newcommand\scf{{\mathscr F}}

\newcommand\scd{{\mathscr D}}

\newcommand\sch{{\mathscr H}}

\newcommand\gra{\alpha}
\newcommand\grb{\beta}

\newcommand\grk{\kappa}

\newcommand\vphi{\varphi}
\newcommand\grr{\rho}
\newcommand\grs{\sigma}

\newcommand\grz{\zeta}

\newcommand\rat{{\mathbb Q}}
\newcommand\reals{{\mathbb R}}
\newcommand\comp{{\mathbb C}}
\newcommand\zed{{\mathbb Z}}

\newcommand\pn[1]{{\mathbb P}^{#1}}
\newcommand\proj[1]{{\mathbb P}({#1})}
\newcommand\Proj{{\mathbb P}}

\newtheorem{theorem}{Theorem}[section]
\newtheorem{lemma}[theorem]{Lemma}
\newtheorem{corollary}[theorem]{Corollary}
\newtheorem{conjecture}[theorem]{Conjecture}
\newtheorem{clam}[theorem]{Claim}
\newtheorem{prop}[theorem]{Proposition}

\theoremstyle{definition}
\newtheorem{definition}[theorem]{Definition}

\newtheorem{defre}[theorem]{Definition-Remark}
\newtheorem{pargrph}[theorem]{}
\newtheorem{examp}[theorem]{Example}
\newtheorem{examps}[theorem]{Examples}

\newtheorem{question}[theorem]{Question}

\newenvironment{example*}{\begin{examp}}{\end{examp}}
\newenvironment{examples*}{\begin{examps}}{\end{examps}}
\newenvironment{definition*}{\begin{definition}}{\end{definition}}
\newenvironment{question*}{\begin{question}}{\end{question}}

\newenvironment{prgrph*}[1]{\indent\begin{pargrph}{\bf #1.}}{\end{pargrph}}
\newenvironment{defre*}{\begin{defre}}{\end{defre}}
\newenvironment{MM*}{\begin{MM}}{\end{MM}}

\theoremstyle{remark}
\newtheorem{re}[theorem]{Remark}

\newenvironment{rem*}{\begin{re}}{\end{re}}

\begin{document}

\title{A view on extending  morphisms from ample divisors}

 
\author[M.C. Beltrametti]{Mauro C. Beltrametti}
\address{Dipartimento di Matematica\\
Via Dodecaneso 35\\
I-16146 Genova, Italy}
\email{beltrame@dima.unige.it}
\author[P. Ionescu]{Paltin Ionescu}
\thanks{The second named author was supported by the Italian Programme ``Incentivazione alla mobilit\` a
 di studiosi stranieri e italiani residenti all'estero''.}
\address{University of Bucharest, Faculty of
Mathematics and Computer Science\\
14 Academiei str.\\
RO--010014 Bucharest\\
and 
Institute of Mathematics of the Romanian Academy\\
P.O. Box 1--764, 
RO 014700 Bucharest\\
Romania
}
\email{Paltin.Ionescu@imar.ro}
\dedicatory{Dedicated to Andrew J. Sommese on his 60th birthday}

\maketitle

\begin{abstract} The philosophy  that  ``a projective manifold is more special than any of its smooth hyperplane sections" was one of  the classical principles of projective geometry.
Lefschetz type results and related vanishing theorems were among the typically used techniques.
We shall survey most of the problems, results and conjectures in this area, using the modern setting of ample divisors, and (some aspects of) Mori theory.

\end{abstract}

\tableofcontents

\section{Introduction}\label{Intro}\addtocounter{subsection}{1}\setcounter{theorem}{0}

In the context of classical algebraic geometry, consider a given embedded complex projective manifold $X\subset \pn N$. One of the typically used techniques was to replace $X$ by some of its smooth hyperplane sections, $Y\subset X$. Thus, the dimension of $X$ is decreased and classification results may be obtained inductively. The efficiency of the method depends on the possibility of transferring some known special properties from $Y$ to $X$. In general, given $Y \subset \pn{N-1}$, there is no \textit{smooth} $X\subset \pn N$ such that $Y$ is one of its hyperplane sections. One can say that \textit{$X$ is more special then $Y$}.

The present paper is a survey of contemporary aspects of the hyperplane section technique. 
A first important ``modern'' incarnation of the above principle is given by Lefschetz's theorem, showing that the topology of $Y$ 
strongly reflects that of $X$ (see \cite{AndrF}). From a geometrical point of view, we are usually given some regular map, say 
$p: Y\to Z$, making $Y$ \textit{special}; e.g., a Fano fibration. We would like to extend this map to $X$. 
It was discovered by Sommese, in his innovative early paper \cite{So1}, that the extension is always possible, if one only assumes that the general fiber of $p$ has dimension 
at least two. His proof is based on Lefschetz's theorem and on (very much related) vanishing results of Kodaira type. In the same paper, Sommese 
showed that when $p$ is smooth and extends, the dimension of $Z$ cannot be too large. It soon became clear that the extension problem is much
harder when $\dim Y\leq \dim Z+1$. Fujita \cite{FuLef} further refined some of the techniques and considered new applications e.g., when $p$ is a 
$\pn d$-bundle or a blowing-up. In the case of three folds, fine results were found by B\u adescu   \cite{BadAmp1, BadAmp2, BadAmp3},
when $Y$ is a $\pn 1$-bundle over a curve and by Sommese \cite{Somin, SoVar}, when $Y$ is not relatively minimal.
It is worth pointing out that the classical context of hyperplane sections was gradually replaced by the more general situation when $Y\subset X$ is 
merely an ample divisor, and no projective embedding of $X$ is given. This is a substantial generalization, since in the new setting the normal 
bundle of $Y$ in $X$ is \textit{not} specified.

The appearance of Mori theory made possible a change of both the point of view and the techniques (see \cite{Mo1,KMM}). The isomorphism between the Picard 
groups of $X$ and $Y$ given by the Lefschetz theorem leads to an inclusion between the Kleiman--Mori cones $\overline{NE}(Y)$ and 
$\overline{NE}(X)$. As is well known, faces of these cones describe non-trivial morphisms defined on $Y$ and $X$, respectively. So, the original 
question of extending maps from $Y$ to $X$ translates into a comparison problem between these cones. Ideally, when the two cones are equal, all
morphisms from $Y$ extend to $X$ (see e.g., \cite{Ko, Wcontr,BFS1, anoc} and Section 8 for results in this direction, usually when $X$ or $Y$ are Fano manifolds).
In the general case, what we can hope for is to extend the contraction of an extremal ray of $Y$ (cf. \cite{IoNef, Occ}).
This is not always possible, but very few counterexamples are known (see Section 4 for this intriguing aspect). The techniques used 
in this setting are the cone theorem, due to Mori, and the contraction theorem, due to Kawamata--Reid--Shokurov, combined with the well behaved
deformation theory of families of rational curves \cite{KoBook, Debarre, KMM}. See also \cite{An,IoGen, Ko} and \cite{Wcontr} for some useful facts about special families of rational curves, coming from extremal rays.

General results on extending morphisms are discussed in Section 3; in Section 5 we concentrate on the special situation when $p$ is a $\pn
d$-bundle or a blowing-up. We pay special attention to the case of 
$\pn 1$-bundles, which is the most difficult. A (still open) main conjecture on the subject is stated and various related facts are proved in Section 7. The afore mentioned results by B\u adescu and Sommese on three folds are recovered in Section 6, using the Mori theory point of view (cf.\ \cite{IoGen}). In the last section we discuss the ascent of some good properties from $Y$ to $X$: e.g., being uniruled, or rationally connected, or rational, etc.

We have tried to write a complete and coherent exposition, also accessible to the nonspecialist. We included several new proofs and sometimes 
substantial simplifications of the original arguments. Several possible generalizations are mentioned at the end of the paper, together
 with appropriate references to the existing literature.

\section{Background material}\label{Background}

We work over the
complex field $\comp$. Throughout the paper we deal with irreducible, reduced, projective varieties $X$. 
We use the term manifold if $X$ is moreover assumed to be smooth.
We denote by $\sO_X$ the structure sheaf of $X$. For any coherent sheaf $\sF$ on $X$, $h^i(\sF)$ denotes the complex dimension of $H^i(X,\sF)$. If $p: X\to Y$ is a morphism, we write $p_{(i)}$ for its $i$-th direct image.

Let $L$ be a line bundle on $X$.  $L$ is said to be {\em
numerically effective} ({\em nef}, for short) if $L\cdot C\geq 0$ for all
effective curves $C$ on $X$. We say that $L$ is {\em strictly nef} (or {\em numerically positive}) if  $L\cdot C> 0$ for all effective curves $C$ on $X$. $L$ is said to be {\em big} if $\grk(L)=\dim
X$, where $\grk(L)$ denotes the Iitaka dimension of $L$. If $L$ is nef
then this is equivalent to $c_1(L)^n>0$, where $c_1(L)$ is the first Chern
class of $L$ and $n=\dim X$. The pull-back $\iota^*L$ of a line bundle $L$ on $X$ by an embedding
$\iota:Y\hookrightarrow X$ is denoted by $L_Y$. We denote by $N_{Y/X}$ the normal bundle of $Y$ in $X$ and by $K_X$ the canonical bundle of a smooth  variety $X$.
 
We use standard notation 
from algebraic geometry, among which we recall the following ones:
\begin{enumerate}
\item[$\bullet$] $\approx$, the linear
equivalence of line bundles; $\sim$, the numerical
equivalence of line bundles;
\item[$\bullet$] $|L|$, the complete linear system associated to a line bundle $L$;
\item[$\bullet$] $\grk(D)$, the Iitaka dimension of the line bundle associated to a
$\mathbb Q$-Cartier divisor $D$ on $X$; and
$\grk(X):=\grk(K_X)$, the Kodaira dimension of $X$, for $X$
smooth.
\item[$\bullet$] $\pi_i(X)$, the $i$-th homotopy group, omitting the base point when its choice is irrelevant.
\end{enumerate}

 $\Proj^n$ denotes the projective $n$-space, $\sQ^n\subset \Proj^{n+1}$  denotes  the $n$-dimensional hyperquadric. For a vector bundle $\sE$, we write $\Proj(\sE)$  for the associated projective bundle and $\xi _\Proj$, or $\xi_X$ when $X=\Proj(\sE)$, for the tautological line bundle, using the Grothendieck convention.   

Line bundles and divisors are used with little (or no)
distinction. We almost always use the additive notation.
We say
that a line bundle $L$ is {\em spanned} if it is spanned,
i.e., globally generated, at all points of $X$ by $H^0(X, L)$.

\subsection{Setting up and motivation}\label{motivations} Let $X$ be a projective manifold and let $Y\subset X$ be a smooth ample divisor. It is a natural classical question  to try  to understand how the structure of $Y$ determines the one of $X$.

More precisely, given a surjective morphism $p=p_{|D|}:Y\to Z$ associated to a linear system $|D|$,  we look  for  a linear system $|\overline{D}|$ on $X$ defining a regular map $\overline{p}=p_{|\overline{D}|}:X\to W$ onto a projective variety $W$, such that the following diagram
\begin{equation}\label{diagram}
\xymatrix{
Y \ar[d]_p \ar @{^{(}->}[r] & X \ar[d]^{\overline{p}} \\
Z \ar[r]_\gra & W }
\end{equation}
commutes.  If  the morphism $\gra:Z\to \gra(Z)$ is finite  we say that  $\overline{p}$ is a {\em lifting} of $p$. If $\overline{p}_{|Y}=p$, that is if $\gra:Z\to \gra(Z)\subset W$ is an isomorphism onto its image, we say that $\overline{p}$ is a {\em strict lifting} of $p$, or that $p$ is {\em extendable} to $\overline{p}$. Note that this is always the case   whenever  the restriction map $H^0(X, \overline{D})\to H^0(Y, D)$ is surjective.  Note also that this further condition will be a posteriori satisfied  in our setting (see the proof of  Theorem~\ref{six}).

Assume that the morphism $p$ has a lifting $\overline{p}$. Up to taking the Remmert--Stein factorization, we can always  assume that $p$ has connected fibers and $Z$ is normal. Therefore, by using the ampleness of $Y$ in $X$, it is a standard fact that one of the following holds:
\begin{enumerate}
\item[(1)]  $\dim Y-\dim Z\geq 1$ and  $\gra:Z\xrightarrow{\sim} W$ (in particular $p$ is extendable);
\item[(2)] $p$, $\overline{p}$ and  $\gra:Z\rightarrow \gra(Z)$ are birational; so, $\gra$ is the normalization morphism;
\item[(3)] $p$ is birational and $\dim X -\dim  W =1$; in this case $\gra :Z\to \gra(Z)$ may be of degree $\geq 2$. 
\end{enumerate}

A simple example is obtained as follows. Consider $X:=\pn 1\times \pn {n-1}$ embedded in $\pn N$ by $\sO(2,1)$, $n\geq 4$.
By Bertini's theorem, we can choose a hyperplane $H$ in $\pn N$  such that the restriction $Y$ of $H$ to $X$ is a smooth ample divisor. Then we get a commutative diagram
\[\xymatrix{
Y \ar[d]_p \ar @{^{(}->}[r]  \ar[dr] & X=\pn 1\times \pn {n-1} \ar[d]^{\overline{p}} \\
Z \ar[r]_\gra &\pn {n-1} }
\]
where $p$ and $\gra$ are given by the  Remmert--Stein factorization of the restriction $\overline{p}_{|Y}:Y\to \Proj^{n-1}$ and $\overline p$ is the 
natural projection.
Note that the morphism $\gra$ is finite of degree two.
Moreover, $p$ is not an isomorphism. Indeed, assume otherwise. Then $Y\to  \pn {n-1}$ is a two-to-one finite covering, so that  it induces an isomorphism $\Pic Y\cong \zed$ (see \cite[II, 7.1.20]{Laz} for details and complete references). On the other hand, $\Pic Y\cong \Pic X\cong \zed\oplus\zed$ by the Lefschetz theorem; a contradiction.

\smallskip

If the morphism $p$ is extendable, our aim is to describe $X$ by using the structure morphism $\overline{p}$. The occurrence that $p$ is not extendable  forces  $X$ to satisfy geometric constraints which, in turn, make $X$ special enough to be completely classified.


As a typical example, consider the following natural question, formulated in the classical context: 
\begin{question}\label{zeroQ} Let $X$ be an $n$-dimensional manifold embedded in a projective space $\pn N$. Assume that a smooth hyperplane section, $Y=X\cap H$, of $X$ is a  $\pn d$-bundle  over  some manifold $Z$, such that the fibers are linearly embedded. Does it  follow that the bundle projection $p:Y\to Z$ extends to $X$ giving a $\pn {d+1}$-bundle projection $\overline{p}:X\to Z$?\end{question}

As soon as $n\geq 4$, the (positive) answer to this question  relies on some non-trivial results from the deformation theory of rational curves. 
It turns out that the  key-fact is the  condition $H\cdot  f=1$, where $f$ is a line in a fiber $\pn d$ of $p:Y\to Z$, i.e., $f$ is a linear $\pn 1$ with respect to the embedding of $X$ in 
$\pn N$ given by $H$. Moreover, the above can happen only if $2\dim Z\leq \dim X$ (Proposition~\ref{Qnot}).



\subsection{Special varieties}\label{specialvar} Let $X$ be a projective manifold of
dimension $n$. We say that $X$ is a {\em Fano manifold} 
if  $-K_X$ is ample; its {\em index}, $i$, is the largest positive integer such that $K_X\approx -i\sL$ for some ample line bundle $\sL$
on $X$. Let $L$ be a given ample line bundle on $X$. We say that $(X,L)$ is a  {\em del Pezzo  variety}
(respectively a {\em Mukai variety}) in the adjunction theoretic sense if $K_X \approx - (n - 1)L$ (respectively
$K_X \approx - (n -2)L)$.
Note that del Pezzo manifolds are completely described by Fujita \cite[I, Section 8]{FuBook}.
We refer to Mukai \cite{Mu1} for results on Mukai varieties.


We say that $(X,L)$ is a {\em scroll} 
over a normal variety $Z$ of dimension $m$ if there exists a surjective morphism with connected
fibers $p : X \to Z$, such that $K_X + (n - m + 1)L \approx p^{\ast} {\sL}$
for some {\em ample} line bundle $\sL$ on $Z$.

We refer to \cite{BSCe} and \cite[Sections 14.1, 14.2]{Book} for relations between the
adjunction theoretic and the classical definition of scrolls.

Let $X$ be a projective manifold and let $p:X\to Z$ be a surjective morphism onto another manifold, $Z$. We say that $X$, $p:X\to Z$ is a $\pn d$-{\em bundle} if each closed fiber of $p$ is isomorphic to the projective space $\pn d$. We also say that  $X$, $p:X\to Z$ is a {\em linear $\pn d$-bundle} if  $X={\mathbb P}(\sE)$ for some rank $d+1$  vector bundle $\sE$ on $Z$.


We  say that  $X$, $p:X\to Z$ is a {\em conic fibration} over a normal projective variety
$Z$ if every fiber of the morphism $p$ is a conic, i.e., it is isomorphic  to the zero scheme of a
non-trivial section of
$\sO_{\pn 2}(2)$. Note that the above definition is equivalent to saying that there exists a rank $3$ vector
bundle $\sE$ over $Z$ such that its projectivization $\widetilde{p}:{\mathbb P}(\sE)\to Z$ contains $X$
embedded over $Z$ as a divisor whose restriction to any fiber of $\widetilde{p}$ is an element of $|\sO_{\pn 2}(2)|$. The push-forward $p_*(-K_X)$ can be taken as the above $\sE$. It is a standard fact to show that  $p:X\to Z$ is a flat morphism; since $X$ is smooth, it follows that the base $Z$ is smooth, too. 



\subsection{Lefschetz-type and vanishing results}\label{LefVan} 
A basic tool  for dealing with the  problems discussed above  are
 Lefschetz's theorems, which, in turn, are very much related (in fact, almost equivalent) to  vanishing results of Kodaira type  
 (see \cite[I, Chapters 3, 4]{Laz} for a nice general presentation  and complete references). See \cite{AndrF} for the classical statement of Lefschetz's theorem.

\begin{theorem}{\em(Hamm--Lefschetz theorem)}\label{LefThm}  Let
$L$ be an ample line bundle on a projective
manifold, $X$, and let $D \in |L|$.
Then given any
point $x\in D$ it follows that the $j$-{\em th} relative homotopy group,
$\pi_j(X,D,x)$, vanishes for $j\le \dim X - 1$.   In particular, the restriction mapping, $H^j(X,\zed) \to H^j(D,\zed)$ is an isomorphism for
$j\le \dim X - 2$,  and is injective with
torsion free cokernel for $j=\dim X -1$.
 \end{theorem}

\begin{theorem}{\em(Barth--Lefschetz theorem)}\label{BarthLefThm} Let   $Y$ be a
  connected  submanifold of a projective manifold, $X$. Let $n=\dim X$,
$m=\dim Y$. Assume that  $N_{Y/X}$  is ample. Then  for any $x\in
Y$, we have $\pi_j(X,Y,x)=0$ for $j\leq 2m-n+1$.
 In particular, under the natural map we have $\pi_1(Y,x)\cong   \pi_1(X,x)$ if $2m-n\geq 1$. Moreover:
\begin{enumerate} 
\item[{\rm (i)}]  If \ $2m-n=1$, the restriction map 
$ r:\Pic X\to\Pic Y$
is injective with torsion free cokernel; and 
\item[{\rm (ii)}]  If \ $2m-n\geq 2$, then $\Pic X \cong \Pic Y$ via $r$.
\end{enumerate}\end{theorem}


Kawamata and Viehweg showed that the Kodaira vanishing theorem holds for any nef and big line bundle (see e.g., \cite[Sections 1--2]{KMM}).

\begin{theorem}{\em(Kawamata--Viehweg vanishing theorem)}\label{KVvan} Let $X$ be a projective manifold of dimension $n$, and let $D$ be a nef and big divisor on $X$. Then
\[H^i(X,\sO_X(K_X+D))=0\quad  \mbox{for }i>0.\]
\end{theorem}

\subsection{Basic facts from Mori theory}\label{MT} Let us recall some definitions and a few facts from Mori theory we need. Basic 
references for details are \cite{Mo1,Mo2}, and \cite{KMM}.
Let $X$ be a connected normal projective variety of dimension $n(\geq 2)$.
\begin{trivlist} \leftskip = 20 pt \hangindent = 20 pt \hangafter =1
\item[$\bullet$] ${\rm Num}(X) = \Pic X /\sim$;
\item[$\bullet$] $N^1(X) = {\rm Num}(X) \otimes \reals$;
\item[$\bullet$] $N_1(X) = (\{1$-{\rm cycles}\}$/\sim) \otimes \reals$;
\item[$\bullet$] $NE(X)$, the convex cone in $N_1(X)$ generated by the effective 1-cycles;
\item[$\bullet$] $\overline{NE} (X)$, the closure of $NE(X)$ in $N_1(X)$ with respect to the Euclidean
topology;
\item[$\bullet$] $\varrho (X) = \dim_\reals N_1(X)$, the {\em Picard number} of $X$;
\item[$\bullet$] $\overline{NE}(X)_{D\geq 0} = \{ \grz \in \overline{NE}(X)\mid  \grz \cdot D \ge 0\}$ for given $ D \in \Pic X
\otimes \rat$;
\item[$\bullet$] $\overline{NE}(X)_{D<0} = \{ \grz \in \overline{NE}(X)\mid  \grz \cdot D < 0\}$ for given $ D \in \Pic X
\otimes \rat$;
\item[$\bullet$] $\overline{NE}(X)_{D\leq 0} = \{ \grz \in \overline{NE}(X)\mid  \grz \cdot D \leq 0\}$ for given $ D \in \Pic X
\otimes \rat$;
 
\item[$\bullet$] ${\rm Nef}(X)$, the dual cone of  $\overline{NE}(X)$, namely, the cone in $N^1(X)$ spanned by classes of nef divisors. 
\end{trivlist}

If $\gamma$ is a $1$-dimensional cycle in $X$ we denote by  $[\gamma]$ its
class in $\overline{NE}(X)$. Note that the vector spaces $N^1(X)$ and $N_1(X)$ are dual to each other via the usual
intersection of cycles ``$\,\cdot\,$".

Assume that $X$ is smooth.  We  say that a half line $R = \reals_+[\grz]$ in $\overline{NE} (X)$, where $\reals_+ =\{x \in \reals\mid x > 0\}$, is an {\em extremal  ray} if  $K_X \cdot \grz < 0$ and
$\grz_1, \grz_2 \in R$ for every $\grz_1, \grz_2 \in \overline{NE}(X)$ such that $\grz_1 + \grz_2 \in R.$

An extremal ray $R = \reals_+[\grz]$ is {\em nef} if $D \cdot \grz \ge 0$ for
every effective divisor $D$ on $X$. An extremal ray which is not nef is said to
be {\em non-nef}.

Let $D \in \Pic X \otimes \rat$ be a nef $\rat$-divisor, $D\not \sim 0$. Let
\[F_D : = D^{\perp} \cap (\overline{NE}(X) \setminus \{0\}),\]
where ``$\perp$" means the orthogonal complement of $D$ in $N_1(X)$. Then $F_D$
is called a {\em good extremal face} of $\overline{NE}(X)$ and
$D$ is the {\em supporting hyperplane} of $F_D$, if $F_D$ is entirely contained in the set
 $\{\grz \in N_1 (X)\mid K_X \cdot \grz < 0\}$. An extremal ray is
a $1$-dimensional good extremal face.  Indeed, for any extremal ray $R$ there exists a
nef $D \in \Pic X \otimes \rat$ such that $R = D^{\perp} \cap
(\overline{NE}(X) \setminus \{0\})$. 

\begin{theorem}{\rm (Mori cone theorem)}\label{ConeTheorem} Let $X$ be a projective manifold of dimension $n$. Then there exists a countable set of curves
$C_i, i \in I$, with $K_X \cdot C_i < 0$, such that one has the decomposition
\[\overline{NE} (X) = \sum_{i \in I} \reals_+ [C_{i}] + \overline{NE}(X)_{K_X\geq 0}.\]
The decomposition has the properties: 
\begin{enumerate}
\item[{\rm (i)}]  the set of  curves $C_{i}$ is minimal, no smaller set is sufficient to generate
the cone;
\item[{\rm (ii)}]  given any neighborhood $U$ of $\overline{NE}(X)_{K_X\geq 0}$, only
finitely many $[C_{i}]$'s do not belong to $U$.
\end{enumerate}
The semi-lines $\reals_+ [C_{i}]$ are the extremal rays of $X$. Moreover,  the
curves $C_{i}$ are {\rm (}possibly singular{\rm )} reduced irreducible rational curves which satisfy the
condition $1 \le -K_X \cdot C_i \le n + 1$. 
\end{theorem}

\begin{theorem}{\rm (Kawamata--Reid--Shokurov base point free theorem)}\label{KRS} Let $X$ 
be a projective manifold of dimension $n \ge 2$. Let $D$ be a nef Cartier divisor such that $aD-K_X$ is
nef and big for some positive integer $a$. Then $|mD|$ has no base points  for $m \gg 0$. 
\end{theorem}

It is a standard fact that, for a good extremal face $F_D$, the line bundle $mD-K_X$ is ample for $m\gg 0$. Therefore,
by  Theorem \ref{KRS}, the linear system $|mD|$ is base point free for $m\gg 0$, so
 that it defines a morphism, say $\vphi : X \to W$. By taking $m$ big enough, we may further assume that $W$ is normal and the fibers of $\vphi$ are connected.  Note that
$\varphi_{\ast} \sO_X \cong \sO_W$, the pair $(W, \varphi)$ is unique up to
isomorphism and $D \in \varphi^{\ast}$ $\Pic W$. If $C$ is an irreducible curve on $X$,
then $[C] \in F_D$ if and only if $D \cdot C = 0$, which means $\dim
\varphi (C) = 0$, i.e.,  $\varphi$ contracts the good extremal face $F_D$. We will call  such a
contraction, $\varphi$, the {\em contraction of} $F_D$. If $F_D=R$, $R$ an extremal ray, we will denote ${\rm cont}_R:X\to W$ the contraction 
morphism. Let
\[E := \{x \in X\mid {\rm cont}_R \mbox{ is   not   an isomorphism  at } x\}.\]
Note that $E$ is the locus of curves whose numerical class is in $R$. We will refer to $E$
simply as the {\em locus} of $R$. 

If $X$ is smooth we define the {\em length of an extremal ray}, \[\length(R) = \min \{-K_X \cdot C\mid C  
\mbox{ rational  curve}, [C] \in R \}.\] 
Note that the cone theorem yields the
bound $0<\length(R)\le n+1$.  We will also use the notation $\length(R)=\ell(R)$. We say that a rational curve $C$ generating an extremal ray $R=\reals_+[C]$ is a {\em minimal curve} if $\ell(R)=-K_X\cdot C$.

The following useful  inequality is inspired by Mori's bend and break 
(cf.\ \cite[Theorem~0.4]{IoGen}, and also  \cite[Theorem~1.1]{Wcontr}, \cite[Corollary~IV.2.6]{KoBook}).

\begin{theorem}\label{ExRayIW}  Let $X$ be a projective manifold of dimension
$n$. Assume that $K_X$ is not nef and let $R$ be an extremal ray on $X$ of length $\ell(R)$. Let $\grr$ be the
contraction of $R$ and let $E$ be any irreducible component of the locus of $R$. Let $\Delta$ be any irreducible
component of any fiber of the restriction, $\grr_E$, of $\grr$ to $E$. Then
\[\dim E+\dim\Delta\ge  n+\ell(R)-1.\]\end{theorem}

By combining the theorem above with a result due to Ando \cite{An}, Wi\'sniewski  \cite[Theorem~(1.2)]{Wcontr} showed the following, which  
plays an important role in the sequel.

\begin{theorem}\label{AW} {\rm (Ando \cite{An}, Wi\'sniewski \cite{Wcontr})} Let $X$ be a projective manifold of dimension $n\geq 3$. Assume that $K_X$ is not nef. Let $\vphi:X\to Z$ be the contraction morphism of an
extremal ray $R$. If every fiber of $\vphi$ has dimension at most one, then $Z$ is smooth and either
$\vphi$ is the blowing-up of a smooth codimension two subvariety of $Z$, or $\vphi$ is a conic
fibration.\end{theorem}

Ando \cite[(3.10), (2.3)]{An} proved the theorem above assuming that the locus $E$ of $R$ satisfies the
condition  $\dim E\geq n-1$. From  the inequality  of Theorem~\ref{ExRayIW} it follows that this is
the case. Indeed, if $\dim E\leq n-2$, for any irreducible component $\Delta$ of any fiber of the restriction of $\vphi$ to $E$, we would have
\[\dim \Delta \geq \ell(R)+1\geq 2,\] contradicting   the fibers  dimension assumption.


\subsection{Families of rational curves}\label{families}
We follow the  notation in \cite{KoBook}, to which we refer for details; see also \cite{Debarre}. Let $X$ be a projective manifold. By ${\rm Hom}_{\rm bir}(\pn 1,X)$ we denote the scheme parameterizing morphisms from $\pn 1$ to $X$ which are birational onto their image. We will denote by $[f]$ the point of  ${\rm Hom}_{\rm bir}(\pn 1,X)$ determined  by such a morphism $f:\pn 1\to X$.


A reduced, irreducible subvariety $V\subset {\rm Hom}_{\rm bir}(\pn 1,X)$ determines a {\em  family of  rational curves} on $X$. We let $\scf$ be the universal family, restricted to $V$, with $p: \scf\to V$ and $q:\scf\to X$ the natural projections. We call the image of $q$ the {\em locus of the family}, denoted by ${\rm Locus}(V)$. A {\em covering family} is a family satisfying ${\rm Locus}(V)=X$.
We say that a family $V$, closed under the action of  ${\rm Aut}(\pn 1)$, is {\em unsplit} if  the image of $V$ in ${\rm  Chow}(X)$ under  the natural morphism $[f]\to [f(\pn 1)]$ is closed. 
In general, the closure of the image of $V$ in ${\rm Chow}(X)$ determines a {\em family of rational 1-cycles} on $X$.
If $x\in X$ is a fixed (closed) point, we denote by $V_x$ the closed subfamily of $V$ consisting of morphisms sending a fixed point $O\in \pn 1$ to $x$.
We say that $V$ is {\em locally unsplit} if, for $x\in {\rm Locus}(V)$ a general point, the family $V_x$ is unsplit.
A family $V$ of rational 1-cycles on $X$ is {\em quasi-unsplit} if any two irreducible components of cycles in $V$ are numerically proportional. Such families typically arise from cycles belonging to an extremal ray.



\section{General results}\label{general}\addtocounter{subsection}{1}\setcounter{theorem}{0}

We discuss throughout  this  section some general results  on   extending morphisms  $p:Y\to Z$ from ample (smooth)  divisors $Y$ of a manifold $X$.

To begin with, let us prove two early theorems due to  Sommese \cite{So1} (see also \cite[(5.2.1), (5.2.5)]{Book}) that marked the starting point of the subject.  The first one  shows that  the morphism $p$ is always extendable   whenever $\dim Y-\dim Z\geq 2$. The second one 
gives the restriction that $\dim X\geq 2\dim Z$  for a smooth $p:Y\to Z$  to extend.

\begin{theorem}\label{ExtensionThm} {\rm (Sommese \cite{So1})} Let $Y$ be a smooth  ample  divisor on
a projective manifold  $X$. Let
$p: Y\to Z$ be a surjective morphism.  If \[\dim Y-\dim Z \ge 2,\]
then $p$ extends to a surjective  morphism $\overline{p} : X \to Z$.
 \end{theorem}
\proof 
Let $\dim X =: n$. Without loss of generality it can be assumed that $\dim Z \ge 1$. Thus
 we have that $n \ge 4$ and therefore by  Lefschetz theorem
 we see that the restriction map gives an  isomorphism,
$\Pic X\cong \Pic Y$. Moreover, by Remmert--Stein factorizing $p$ it can be
assumed  that $Z$ is normal and $p$ has connected fibers.

Let $\sL$ be a very ample line bundle on $Z$. Since $\Pic X\cong \Pic Y$
there exists an  $H\in \Pic X$ whose  restriction,
$H_Y$,  to $Y$ is isomorphic to $p^*\sL$.
Let $L:=\sO_X(Y)$. We claim that
\begin{equation}\label{vanishing} H^1(Y,(H-tL)_Y)=0 \quad  \mbox{for } t\geq 1.\end{equation}
By Serre duality we are
reduced to showing $H^{n-2}(Y,K_Y+ (tL-H)_Y)=0$ for
$t\ge 1$.  

From the relative form of the Kodaira vanishing theorem (see e.g., \cite{KMM,EV,SS}) we see that
\[p_{(j)}(K_Y+ (tL-H)_Y) = p_{(j)}(K_Y+t L_Y)\otimes (-\sL)=0\]
 for  $j\geq 1$. 
Using the Leray spectral
sequence we deduce that 

\[H^{n-2}(Y, K_Y+(tL-H)_Y)=H^{n-2}(Z, p_*(K_Y+(tL-H)_Y)).\]
The last group is zero, since $n-2=\dim Y-1> \dim Z$ by our assumption.
This shows  (\ref{vanishing}).

Now consider the exact sequence
\[0\to K_X\otimes (L-H)\otimes (tL)\to K_X\otimes (L-H)\otimes (t+1)L\to K_Y\otimes(L-H)_Y\otimes (tL_Y)\to
0.\]
By (\ref{vanishing}) we have $H^1(Y, (H-L)_Y-tL_Y)=0$ for $t\geq 0$. Therefore, by Serre duality,
$H^{n-2}(Y, K_Y\otimes(L-H)_Y\otimes (tL_Y))=0$ for $t\geq 0$. Thus the exact sequence above gives an injection
 of $H^{n-1}(X, K_X\otimes (L-H)\otimes (tL))$ into  $H^{n-1}(X, K_X\otimes (L-H)\otimes (t+1)L)$. By Serre's
vanishing theorem, $H^{n-1}(X, K_X\otimes (L-H)\otimes (t+1)L)=0$ for $t\gg 0$. Therefore we conclude that
\[H^{n-1}(X, K_X\otimes (L-H)\otimes (tL))=H^1(X, H-(t+1)L)=0 \quad \mbox{for } t\geq 0.\]
Hence in particular $H^1(X, H-L)=0$, so that we have a surjection
\[H^0(X, H)\to H^0(Y, H_Y)\to 0.\]
 Since $H_Y\cong p^*\sL$ with $\sL$ very ample on $Z$, we infer that there exist
$\dim Z+1$ divisors $D_1,\ldots,D_{\dim Z+1}$ in $|H|$ such that $D_1\cap\cdots\cap D_{\dim
Z+1}\cap Y=\emptyset$. Since $Y$ is ample it thus follows that 
 $\dim(\cap_{i=1}^{\dim Z+1} D_i)\leq 0$.
We claim that $H$ is spanned by its global sections. If $\cap_{i=1}^{\dim Z+1} D_i=\emptyset$, then
${\rm Bs}|H|=\emptyset$. If
$\cap_{i=1}^{\dim Z+1} D_i\neq\emptyset$, then
\[\dim(D_1\cap\cdots\cap D_{\dim Z+1})\geq \dim X-\dim Z-1\geq 2.\]
 This contradicts the above inequality.

Let $\overline{p} :  X \to  {\mathbb P}^N$ be the map associated to $H^0(X, H)$. Ampleness of $Y$ yields that $\overline p(X)=p(Y)$, so $\overline{p}|_Y=p$ and we are done. \qed

\begin{theorem}\label{basere} {\rm (Sommese \cite{So1})} Let $Y$ be a smooth ample divisor on a projective manifold, $X$. Let $p: Y\to Z$ be a morphism of maximal rank onto a normal variety, $Z$.
  If $p$ extends to a morphism $\overline{p} : X \to Z$, then $\dim X\geq 2\dim Z$. 
 \end{theorem}
\proof  We follow the topological argument from  \cite[Proposition V]{So1}.  Let $S$ be the image of the  set of points where $\overline{p}$ is not of maximal rank. By ampleness of $Y$  the set $S$ is  finite. Let $f=p^{-1}(z)$, $F=\overline{p}^{-1}(z)$ be the fibers of $p$, $\overline{p}$ over $z\in Z\setminus S$ respectively. Let $\dim Z=b$ and let $r=\dim Y-\dim Z$.
From standard results in topology we deduce:
$\pi_j(Y,f)\cong \pi_j(Z)$ for all $j$, $\pi_j(X,F)\cong \pi_j(Z)$ for $j\le 2b-2$ and $\pi_{2b-1}(X,F)\to \pi_{2b-1}(Z)$ is onto.
It follows that $\pi_j(Y,f)\to \pi_j(X,F)$ is an isomorphism for $j\le 2b-2$ and it is onto for $j=2b-1$.
From Theorem~\ref{LefThm}, we have that $\pi_j(Y)\to \pi_j(X)$ is an isomorphism if $j<$ dim$Y$ and is onto for $j=$ dim$Y$.
Consider the following commutative diagram with exact rows:
\[\begin{array}{ccccccccccc}
\cdots&\to&\pi_j(f)& \to &\pi_j(Y)&\to& \pi_j(Y,f)&\to &\pi_{j-1}(f)&\to&\cdots\\
&&\downarrow&  &\downarrow&& \downarrow& &\downarrow&&\\
\cdots&\to&\pi_j(F)& \to &\pi_j(X)&\to& \pi_j(X,F)&\to &\pi_{j-1}(F)&\to&\cdots\\
\end{array}\]

Arguing by contradiction, let us assume $\dim X< 2\dim Z$, or $r<b-1$. It follows that $2r+2\le r+b=$ dim$Y$ and $2r+2 < 2b-1$. This, the above, and the five lemma show that $\pi_j(f)\cong\pi_j(F)$ for $j<2r+2$ and
$\pi_{2r+2}(f)\to \pi_{2r+2}(F)$ is onto.

By Whitehead's generalization of Hurewicz's theorem \cite[Theorem~9, p. 399]{Spanier} we get
$H_j(f,\zed)\cong H_j(F,\zed)$ for $j<2r+2$ and a surjection
\[H_{2r+2}(f,\zed)\to H_{2r+2}(F,\zed)\to 0.\]
By noting that $2r+2=2(\dim f+1)=\dim_\reals F$, this leads to the contradiction
$H_{2r+2}(F,\zed)=0$.
\qed

\begin{rem*}\label{equality} In the situation of Theorem \ref{basere}, if  $\dim X=2\dim Z$, a more refined argument based on results of Lanteri and Struppa \cite{LaSt} shows that   the general fiber $(F,\sO_F(Y))$  of $\overline{p}$ is isomorphic to $(\pn {\dim F}, \sO_{\pn {\dim F}}(1))$. We refer to \cite[(5.2.6), (2.3.9)]{Book} for more details and complete references.
\end{rem*}

From now on, \begin{itemize}
 \item {\em we are  reduced to consider the problem  of extending morphisms  $Y\to Z$ from ample  divisors $Y$ of a manifold $X$ in the hardest case when} $\dim Y-\dim Z\leq 1$. \end{itemize}
 In this setting, we shall use (part of) Mori theory; we compare the cone of curves of $Y$ and $X$ and give results on extending contractions of extremal rays (cf. \cite{Ko,Wcontr,IoNef}).

As noted in \cite[Section 3]{anoc}, the following  useful fact holds true. It is an easy consequence of Theorem \ref{AW}
and a  lemma due to Koll\'ar \cite{Ko}.

\begin{prop}\label{KollarProp} Let $X$ be a projective manifold of dimension $\geq 4$. Assume that $K_X$ is
not nef and let $R=\reals_+[C]$ be an extremal ray on $X$. Let $Y$ be a smooth ample divisor on $X$. If
$(K_X+Y)\cdot C\leq 0$, then $R\subset \overline{NE}(Y)$.
\end{prop} 
\proof By the Lefschetz theorem, 
the embedding $i:Y\hookrightarrow X$ gives an isomorphism $N_1(Y)\cong N_1(X)$, under which we get a natural
inclusion $i_*:\overline{NE}(Y)\hookrightarrow \overline{NE}(X)$.

Let $\vphi:X\to Z$ be the contraction of the extremal ray $R$ and let $E$ be the locus where $\vphi$ is not
 an isomorphism, i.e., the locus of curves whose numerical class is in $R$. If there is a fiber $F\subset X$ of
$\vphi$ whose dimension is at least two, then $Y\cap F$ contains a curve $\gamma$ which generates $R$ in
$\overline{NE}(X)$, and hence $R\subset \overline{NE}(Y)$. 
 Thus we can assume
that every fiber of $\vphi$ has dimension at most one, so that Theorem \ref{AW} applies. 
Therefore  we are done after showing that in each case of \ref{AW} the divisor $Y$ contains a fiber
of $\vphi$.

In the birational case, $E$ is a $\pn 1$-bundle over $\vphi(E)$. Let $F\cong\pn 1$ be a fiber of the bundle
projection $E\to \vphi(E)$. Then $-K_X\cdot F=1$, so that $(K_X+Y)\cdot F\leq 0$ and the ampleness of $Y$ give $Y\cdot
F=1$. Therefore Lemma~\ref{Kollar}(i) below leads to the
contradiction $\dim \vphi(E)\leq 1$, so $\dim X\leq 3$.
In the conic fibration case, for any fiber $F$ of $\vphi$, we have $-F\cdot K_X\leq 2$, and hence we get 
$1\leq Y\cdot F\leq 2$. Thus  Lemma~\ref{Kollar}(ii) gives the contradiction $\dim Z\leq 2$.
\qed

Let us point out the following  consequence of Proposition \ref{KollarProp} (cf. Section \ref{Fano}).

\begin{prop}\label{reph} Let $X$ be a projective manifold of dimension $n\geq 4$, let $H$ be an ample
line bundle on $X$, and let $Y$ be an effective smooth divisor  in $|H|$. Assume that $-(K_X+H)$ is nef. Then
$X$ is a Fano manifold and $\overline{NE}(X)\cong \overline{NE}(Y)$.
\end{prop}
\proof Let $D:=-(K_X+H)$, which is nef. Then $-K_X=H+D$ is
ample, so that $X$ is a Fano manifold. Let $R=\reals_+[C]$ be an extremal ray in the polyhedral cone 
$\overline{NE}(X)$. By assumption, $(K_X+H)\cdot
C\leq 0$. Then Proposition \ref{KollarProp}  applies to give that $R$ is contained in $\overline{NE}(Y)$. So
$\overline{NE}(X)=\overline{NE}(Y)$. \qed

\begin{lemma}\label{Kollar} {\rm (Koll\'ar \cite{Ko})} Let $X$ be a projective manifold.
\begin{itemize}
\item[{\rm (i)}] Let $p:X\to Z$ be a $\pn 1$-bundle over a normal projective variety $Z$. Let $Y\subset X$ be a
divisor such that the restriction $p:Y\to Z$ is finite of degree one. If $Y$ is ample then $\dim Z\leq 1$.
\item[{\rm (ii)}]  Let $p:X\to Z$ be a conic fibration over a normal projective variety $Z$. Let $Y\subset X$ be a divisor such that the restriction  $p:Y\to Z$ is finite of degree two {\rm (}or one{\rm )}. If $Y$ is ample then $\dim Z\leq
2$.\end{itemize}
\end{lemma}
\proof (i) Note that $p_*\sO_X(Y)$ is an ample rank $2$ vector bundle since $Y$ is ample. On the other hand,  the section $\sO_X\to \sO_X(Y)$ gives an extension
\[0\to \sO_Z\to p_*\sO_X(Y)\to L\to 0,\] where $L$ is a line bundle. Thus $c_2(p_*\sO_X(Y))=0$, which contradicts ampleness for $\dim Y\geq 2$.

(ii) If  $p:X\to Z$  has only smooth fibers, then, after a finite base change $Z'\to Z$, we get a $\pn 1$-bundle $p':X'\to Z'$. Now, $p'_*\sO_{X'}(Y')$ has rank $3$ and we get an extension
\[0\to \sO_{Z'}\to p'_*\sO_{X'}(Y')\to \sE\to 0,\] where $\sE$ is a rank $2$ vector bundle.  Thus 
$c_3(p'_*\sO_{X'}(Y'))=0$, which contradicts ampleness for $\dim Z\geq 3$.

If  $p:X\to Z$  has singular fibers, then  let $p':X'\to Z'$ be the universal family of lines in the fibers. The pull-back of $Y$ to $X'$ intersects every line once and it is ample. Moreover, $\dim Z'=\dim Z-1$, thus we are done by (i).
\qed

\begin{corollary}\label{KollarCor}  Let $X$ be a projective manifold of dimension $\geq 4$.  Let $Y$ be a smooth ample divisor on $X$.  Assume that $K_Y$ is not nef and let $R_Y$ be an extremal ray on $Y$. Further assume
that there exists a nef divisor $D$ on $X$ such that $D\cdot R_Y=0$ and  such that $D^\perp$ is  contained in  the set $\{ \grz \in \overline{NE}(X)
\mid (K_X+Y)\cdot \grz \le 0\}$.
 Then there exists an extremal ray  $R \subset
\overline{NE}(X)$ which induces $R_Y$. \end{corollary}
\proof Note that   $D^\perp$ is  a locally polyhedral face of  $\overline{NE}(X)$ since $\{ \grz \in \overline{NE}(X) \mid (K_X+Y)\cdot \grz \le 0\}\subset \overline{NE}(X)_{K_X<0}\cup \{0\}$. Let $D^\perp= \langle R_1,\ldots,R_s\rangle$, $R_i$ extremal rays on $X$, $i=1,\ldots,s$. We can assume $R_Y\neq R_i$ for each $i$ since otherwise we are done. Since
$D\cdot R_Y=0$, it follows that, for some $s'\leq s$, 
\[R_Y\subset \langle R_1,\ldots,R_{s'}\rangle \subseteq \{ \grz \in \overline{NE}(X) \mid (K_X+Y)\cdot \grz \le 0\}.\]
Therefore $R_i\subset \overline{NE}(Y)$, $i=1,\ldots,s'$, by Proposition \ref{KollarProp}. Recalling that $R_Y\neq R_i$,
we thus conclude  that $R_Y$ is not  extremal on $Y$, a
contradiction. \qed 

The following numerical invariant was introduced in \cite{IoNef}. Let $X$ be a projective manifold
of dimension $n\geq 3$. Let $Y$ be a smooth ample divisor on $X$. Let  $R$ be  an 
extremal ray on $Y$ and let $H$ be an ample divisor on $X$. Then define
\[\gra_H(R):=\frac{H\cdot C}{\ell(R)}, \quad C  \mbox{ minimal rational curve generating }  R,\]
and
\[\gra_H(Y):=\min\{ \gra_H(R)\; | \; R\; {\rm extremal \; ray \; on} \; Y\}.\]

The following slightly improves the main result of \cite{IoNef}. 

\begin{theorem}\label{six} Let $X$ be a projective manifold
of dimension $n\geq 4$.  Let $Y$ be a smooth ample divisor on $X$. Assume that $K_Y$ is not nef and let  $R$ be  an  extremal ray on $Y$. Let $p:={\rm cont}_R:Y\to Z$ be the contraction of $R$. The following conditions are equivalent:
\begin{enumerate}
\item[{\rm (i)}] There exists an extremal ray $\scr$ on $X$  which induces $R$ on $Y$;
\item[{\rm (ii)}] There exists a nef divisor $\scd$ on $X$ such that $R=\scd_Y^\perp\cap(\overline{NE}(Y)\setminus\{0\})$;
\item[{\rm (iii)}] $p$ is extendable;
\item[{\rm (iv)}] $p$ has a lifting;
\item[{\rm (v)}] There exists an ample line bundle $H$ on $X$ such that $\gra_H(R)=\gra_H(Y)$.
\end{enumerate}
\end{theorem}
\proof  
 (i) $ \Longrightarrow$ (ii) We have $\scr=\scd^\perp\cap (\overline{NE}(X)\setminus\{0\})$  for some nef divisor $\scd$ on $X$. Restricting to $Y$, gives $R=\scd_Y^\perp\cap(\overline{NE}(Y)\setminus\{0\})$.

 (ii) $ \Longrightarrow$ (iii) There exists an extremal ray  
$\scr \subseteq \scd^\perp\cap\overline{NE}(X)_{K_X+Y<0}$. Then by Proposition~\ref{KollarProp}, $\scr\subset \overline{NE}(Y)$, that is, $\scr$ induces $R$ on $Y$. Replacing $\scd$ by some other nef divisor on $X$, we may assume that $\scr=\scd^\perp\cap (\overline{NE}(X)\setminus \{0\})$. It follows that $m\scd -(K_X+Y)$ is ample for $m \gg 0$ and cont$_{\scr}$ is given by $|m\scd|$. By Kodaira vanishing we have $H^1(X, m\scd -Y)=0$, showing that cont$_{\scr}$ extends cont$_R$.
 
 (iii) $\Longrightarrow$ (iv) is obvious.
 
 (iv) $\Longrightarrow$ (i) The morphism  $\overline{p}:X\to W$ which lifts  $p$  is associated to a complete linear system $|L|$, for some nef line bundle $L$ on $X$. Clearly, $R\subset L^\perp$ and $R\cdot (K_X+Y)=R\cdot K_Y<0$. Hence, in particular, $L^\perp\cap \{ \grz \in \overline{NE}(X) \mid (K_X+Y)\cdot \grz \leq 0\}\neq (0)$. Therefore, since $\{ \grz \in \overline{NE}(X) \mid (K_X+Y)\cdot \grz \le 0\}$ is locally polyhedral
in $\overline{NE}(X)$ (see also the proof of \ref{KollarCor}), there exists an extremal ray $\scr$ on $X$ such that $\scr\subset L^\perp\cap\{ \grz \in \overline{NE}(X) \mid (K_X+Y)\cdot \grz \le 0\}$. We also know  by Proposition \ref{KollarProp} that $\scr\subset L^\perp\cap  (\overline{NE}(Y)\setminus\{0\})$. It  thus follows  that  $p$ contracts $\scr$. Since $p$ is the contraction of an extremal ray, $R$,  we conclude that $R=\scr$  in $ \overline{NE}(Y)$. 

(i) $ \Longrightarrow$ (v)  As noted in the proof of (ii)$\Longrightarrow$ (iii), if  $\scr=\scd^\perp\cap (\overline{NE}(X)\setminus\{0\})$, we have that $m\scd -(K_X+Y):=H$ is ample for $m \gg 0$. Nefness of $\scd$ yields for each extremal ray $ R'=\reals_+[C']$ on $Y$, $C'$ minimal curve, $(K_Y+H_Y)\cdot C' \ge 0$. Moreover, since $R\subset \scd^\perp$, we get $1=\gra_H(R)= \frac{H\cdot C}{\ell(R)}
\leq \frac{H\cdot C'}{\ell(R')}=\gra_H(R')$. This means that $\gra_H(R)=\gra_H(Y)$ is minimal.

 Finally,  to show (v) $ \Longrightarrow$ (i), let  $R:=\reals_+[C]$, $C$ minimal rational curve. 
Let $a=(H\cdot C)$ and let $b=\ell(R)$.  Assume that   $\gra_H(R)=\gra_H(Y)$. 
Then $aK_Y+bH_Y$ is nef on $Y$, since otherwise
$(aK_Y+bH_Y)\cdot C'<0$ for some extremal ray $R'=\reals_+[C']$ on $Y$. This would give
$\frac{a}{b}=\frac{(H\cdot C)}{\ell(R)}>\frac{(H\cdot C')}{\ell(R')}$, contradicting the
minimality of $\gra_H(R)$. 

We claim  that $D:=a(K_X+Y)+bH$ is nef on $X$. 
Assuming the contrary, by Mori's cone theorem there exists an extremal ray
$R'=\reals_+[C']$, such that $(a(K_X+Y)+bH)\cdot C'<0$. Let $\rho$ be the contraction of $R'$. If there is a fiber $F$
of $\rho$ of dimension $\geq 2$, then some curve $\gamma\subset F$ is contained  in $Y$  by the ampleness of  $Y$. Since $(a(K_X+Y)+bH)\cdot \gamma<0$, this contradicts
the nefness of the restriction of $a(K_X+Y)+bH$ to $Y$. Thus we conclude that all fibers of $\rho$ are of
dimension $\leq 1$. By Theorem~\ref{AW} we deduce that $X$ is either a conic fibration, or the
blowing-up of a smooth codimension two subvariety. In this latter case one has $K_X\cdot C'=-1$,
which contradicts $(a(K_X+Y)+bH)\cdot C'<0$. In the first case  one has either $K_X\cdot C'=-2$ or 
 $K_X\cdot C'=-1$ according to whether  $C'$ is a conic or a line. If  $K_X\cdot C'=-1$ we have the same numerical contradiction.  If  $K_X\cdot C'=-2$ we get $Y\cdot C'=1$
by using again the above inequality. It thus follows that $X$ is a $\pn 1$-bundle over a smooth  variety $W$ and
$Y$ is a smooth section. By pushing forward the exact sequence
\[0\to \sO_X\to \sO_X(Y)\to \sO_Y(Y)\to 0,\] we get an exact sequence 
\[0\to \sO_W\to \sE\to L\to 0,\] where $\sE$ is an ample  rank $2$ vector bundle such that $X={\mathbb P}(\sE)$ and $L\cong
N_{Y/X}$ is an ample line bundle. By the Kodaira vanishing theorem we have
$H^1(W,-L)=0$, so that the sequence splits. This contradicts the  ampleness of $\sE$. Thus we conclude that $a(K_X+Y)+bH=D$ is nef.

 Now, note that $D\cdot R=0$ and that $D^\perp$ is  strictly contained in $\overline{NE}(X)_{K_X+Y<0}$.
Thus Corollary~\ref{KollarCor} applies to say that there exists an extremal ray $\scr\subset \overline{NE}(X)$ which induces $R$.
\qed

\begin{corollary}\label{newcor}  Let $X$ be a projective manifold
of dimension $n\geq 4$.  Let $Y$ be a smooth ample divisor on $X$. Assume that $K_Y$ is not nef. Then there exists an extremal ray on $Y$ which extends to an extremal ray  on $X$. In particular, when $Y$ has only one extremal ray, it always extends.
\end{corollary}
\proof
For instance, take  any extremal ray $R$  of $Y$ attaining the minimal value of the invariant  $\gra_Y(R)$ (i.e., $\gra_Y(R)=\gra_Y(Y)$).
\qed



The following theorem is a version of a result of Occhetta \cite[Proposition~5]{Occ} (who states it in the
more general case when $Y$ is the zero locus of a section of  an ample vector bundle $\sE$ on $X$ of the expected
dimension $\dim X-{\rm rank}\sE$). His argument contains an unclear critical  point. With notation as in the
proof below, the conclusion in \cite{Occ} uses in an essential way the fact that $D$ is an adjoint divisor,
i.e., $D=aK_X+\sL$ for some  integer $a>0$ and some {\em ample} line bundle $\sL$ on $X$ (implying   that $D^\perp$ is strictly contained in $ \overline{NE}(X)_{K_X<0}$). 
However, we only know that $\sL_Y $ is ample! See also Remark~\ref{R3.12} below.

\begin{theorem}\label{Occhetta} {\rm (cf. Occhetta \cite{Occ})} Let $X$ be a projective manifold of dimension $\geq 4$.  Let $Y$ be a
smooth ample divisor on $X$.  Assume that $K_Y$ is not nef and let $R$ be an extremal ray on $Y$. Let $C \subset Y$ be a rational curve whose numerical class is in $R$. Assume that the deformations of $C$ in $X$ yield a covering and quasi-unsplit family of rational cycles. Then $R$ is an extremal ray of $X$, too. In particular, the conclusion holds if the deformations of $C$ in $Y$ cover $Y$ and $H\cdot C=1$ for some ample divisor $H$ on $X$.\end{theorem}
\proof  We show first the last claim. Let $\nu:\pn 1\to C$  be the normalization of $C$ and let $g:\pn 1\to Y$ and $f:\pn 1\to X$ be the induced morphisms to $Y$ and $X$ respectively. If $C$ yields a covering family of $Y$, its deformations in $X$ cover $X$, too. Indeed, consider
 the tangent bundle sequence
\[0\to T_Y\to T_{X|Y}\to \sO_Y(Y)\to 0.\] By pulling back to $\pn 1$, we get the exact sequence
\[0\to g^*T_Y\to g^*(T_{X|Y})=f^*T_X\to g^*\sO_Y(Y)\to 0.\]
Since both $g^*\sO_Y(Y)$ and $g^*T_Y$ are nef and hence  spanned, we conclude that $f^*T_X$ is spanned, and
therefore that $C$ induces a covering family on $X$, see \cite[II, Section~3, IV, (1.9)]{KoBook}. Moreover, the condition $H\cdot C=1$ ensures that the deformations of $C$ in $X$ yield an unsplit, hence also quasi-unsplit, family.

Let $D_Y\in \Pic Y$ be some (nef) supporting divisor of $R$, i.e., $R=D_Y^{\perp}\cap (\overline{NE}(Y)\setminus\{0\})$. 
By the Lefschetz
theorem, there exists a divisor $D\in \Pic X$ which restricts to $D_Y$. Let $V$ be a covering and quasi-unsplit family of rational cycles on $X$, containing $C$.  

\begin{clam}\label{Wclaim} $D$ is nef.\end{clam}

\proof Assume that $D\cdot \Gamma<0$ for some irreducible curve $\Gamma$ on $X$. Clearly, we can assume that
$\Gamma$ is not contained in $Y$ since the restriction of $D$ to $Y$ is nef. Since $V$ is a covering family and $Y$ is ample, we
can find a curve $B'$ in $V$ parameterizing curves meeting both $\Gamma$ and $Y$. Let $B$ be the normalization of
$B'$. Consider the base-change diagram
\[  \xymatrix{
\widetilde{S} \ar[r] \ar[rd]_{\pi} \ar @/^1.5pc/[rrr]^{\psi} & S \ar[r] \ar[d] & \scf \ar[r]_{q} \ar[d]^p & X \\
 & B \ar[r]  & V & & }
 \]
where $\scf$ is the universal family and $\widetilde{S}$ is a desingularization of $S$, an irreducible component of $p^{-1}(B)$, whose locus contains $\Gamma$. Note that $\widetilde{S}$ is a ruled surface over the curve $B$. 
Let $A:=\psi(\widetilde{S})\cap Y$ be the trace on $Y$ of the image in $X$ of the surface $\widetilde{S}$. Since
$A$ is an ample divisor on $\psi(\widetilde{S})$, there exists at least one irreducible component, say $\sC$, of $A$ which is
not contracted by ${\rm cont}_R:Y\to Z$.  Let $\widetilde{\Gamma}$, $\widetilde{\sC}$ be two irreducible curves on
$\widetilde{S}$ such that $\psi(\widetilde{\Gamma})=\Gamma$, $\psi(\widetilde{\sC})=\sC$. By the above and the hypothesis that $V$ is quasi-unsplit,
$\widetilde{\sC}$ is not a fiber of $\pi:\widetilde{S}\to B$. 

We can write, for some integers $\varepsilon, \delta_i, \varepsilon>0, \widetilde{\sC}\sim\varepsilon C_0+\sum_i \delta_i F_i$,
where $C_0$ is a section of $\pi$ and each $F_i$ is contained in a fiber of $\pi$. We also have
$\widetilde{\Gamma}\sim \gra C_0+\sum_i \grb_i F_i$, for some integers $\gra$, $\grb_i$, $\gra>0$. Thus 
\[\varepsilon\widetilde{\Gamma}\sim \gra\varepsilon C_0+\varepsilon\sum_i\grb_i F_i\sim \gra \widetilde{\sC}-\sum_i(\gra\delta_i -\varepsilon\grb_i)F_i,\] that
is,  in $\Pic
{\widetilde{S}}\otimes\rat$,  one has $\widetilde{\Gamma}\sim a\widetilde{\sC}+\sum_i b_iF_i$, with $a=\frac{\gra}\varepsilon>0$. 

Let $\widetilde{D}:=\psi^*D$. Since $V$ is quasi-unsplit, $\widetilde{D}\cdot F_i=D\cdot R=0$,  so that
\[D\cdot\Gamma=\widetilde{D}\cdot\widetilde{\Gamma}=\widetilde{D}\cdot
\bigg (a\widetilde{\sC}+\sum_i b_iF_i\bigg)=a(\widetilde{D}\cdot\widetilde{\sC})_{\widetilde{S}}=a(D\cdot \sC)\geq 0\] since $D_Y$ is
nef. This shows the claim.

To conclude  we have to show that $R$ is an extremal ray on $X$; see also the proof of  \cite[Theorem~(5.1)]{BdFL}. By Lefschetz's theorem, the embedding $i:Y\hookrightarrow X$ gives a natural inclusion 
$i_*:\overline{NE}(Y)\hookrightarrow \overline{NE}(X)$. Clearly, $\scr:=i_*(R)$ is $K_X$-negative by the adjunction formula.

Since $R$ is an extremal ray of $\overline{NE}(Y)$, by duality
it corresponds to it an extremal face of  maximal dimension $\varrho(Y)-1$ of  the cone of nef divisors ${\rm Nef}(Y)$.
Therefore we can find $\varrho(Y)-1$ good
supporting divisors of $R$ whose numerical classes are
linearly independent in $N^1(Y)$.
By Claim~\ref{Wclaim}, this implies that such good supporting divisors
extend to divisors on $X$ that are nef, trivial on $\scr$, and whose
numerical classes are linearly independent in $N^1(X)$.
Since there are $\varrho(Y)-1$ of them, and $\varrho(Y) = \varrho(X)$
by the isomorphism $N^1(X) \cong N^1(Y)$,
this implies that $\scr$ is an extremal ray of $\overline{NE}(X)$.
\qed

\begin{rem*} \label{R3.12} In \cite[Proposition 5]{Occ}, the author states the result assuming that $R$ is nef. However, the theorem also applies to non-nef extremal rays of $Y$, see Proposition \ref{Palbup2} below. Note that, even when $R=\reals_+[C]$ is nef, in general $C$ does not define a covering family of $Y$. E.g., take cont$_R$ to be a conic fibration, $C$ being a line in a degenerate fiber.
\end{rem*}
\section{Some convex geometry speculations}\label{Convex}
\addtocounter{subsection}{1}\setcounter{theorem}{0}

First, we recall the following  simple observations, due to B\u adescu.

\begin{lemma}\label{Badlemma} {\rm (\cite[Remark 1), p. 170]{BadAmp1})} On the projective line $\pn 1$, consider  a line bundle $\sO_{\pn 1}(a)$, for some integer $a\geq 2$. Write $a=b+c$, with $b,c>0$. Then there exists a surjective map
\[\sO_{\pn 1}(b)\oplus\sO_{\pn 1}(c)\to \sO_{\pn 1}(a)\to 0.\]
\end{lemma}
\proof Consider the global sections 
\[H^0(\pn 1, \sO_{\pn 1}(a))=\langle u^a,u^{a-1}v,\ldots,v^a\rangle\]
as homogeneous polynomials in two variables $u$, $v$. We have  natural inclusions
\[H^0(\pn 1, \sO_{\pn 1}(b))\subset H^0(\pn 1, \sO_{\pn 1}(a)) \quad \mbox{and}\quad
H^0(\pn 1, \sO_{\pn 1}(c))\subset H^0(\pn 1, \sO_{\pn 1}(a)),\] given by multiplication by $u^{a-b}$ and $v^{a-c}$ respectively. Then we get a surjective map
\[H^0(\pn 1,\sO_{\pn 1}(b))\oplus H^0(\pn 1,\sO_{\pn 1}(c))\to H^0(\pn 1,\sO_{\pn 1}(a))\to 0,\] and injections
\[0\to  \sO_{\pn 1}(b)\stackrel{\grb}{\to} \sO_{\pn 1}(a),  \quad 0\to  \sO_{\pn 1}(c)\stackrel{\gamma}{\to}\sO_{\pn 1}(a).\] Thus $\grb\oplus \gamma:\sO_{\pn 1}(b)\oplus\sO_{\pn 1}(c)\to \sO_{\pn 1}(a)$ gives the requested map. The surjectivity  follows from the commutative square
\[
\xymatrix{
H^0(\pn 1,\sO_{\pn 1}(b))\oplus H^0(\pn 1,\sO_{\pn 1}(c)) \ar[d]_{ev_x}\ar[r] & 
H^0(\pn 1,\sO_{\pn 1}(a)) \ar[d]^{ev_x} \\
(\sO_{\pn 1}(b)\oplus\sO_{\pn 1}(c))_x \ar[r] & (\sO_{\pn 1}(a))_x }
\]
where the vertical arrows are the evaluation maps in a given point $x\in \pn 1$ and $ev_x:H^0(\pn 1,\sO_{\pn 1}(a))\to (\sO_{\pn 1}(a))_x$ is onto by spannedness of $\sO_{\pn 1}(a)$.
\qed

\begin{prop}\label{Badprop} {\rm (\cite{BadAmp1})} Given a vector bundle $\sE=\bigoplus_{i=1}^{n-1}\sO_{\pn 1}(a_i)$, with $a_1\geq 2$, $a_1=b+c$, and $b,c>0$, there exists an exact sequence
\begin{equation}\label{lucianex}
0\to\sO_{\pn 1}\to\sF:= \sO_{\pn 1}(b)\oplus\sO_{\pn 1}(c)\oplus\Big( \textstyle{\bigoplus\limits_{i=2}^{n-1}}\sO_{\pn 1}(a_i)\Big)\to \sE\to 0.\end{equation}
\end{prop}
\proof Lemma \ref{Badlemma} yields a surjective map
\[\sO_{\pn 1}(b)\oplus\sO_{\pn 1}(c)\oplus \Big( \textstyle{\bigoplus\limits_{i=2}^{n-1}}\sO_{\pn 1}(a_i)\Big)\to \sE\to 0,\] whose kernel is the trivial bundle since $\det(\sE)=\det(\sF)$.
\qed

\begin{rem*}\label{lurem} Note that Proposition \ref{Badprop} gives rise to a method to construct ample divisors which are projective bundles over  $\pn 1$. Indeed,  let $Y:={\mathbb P}(\sE)$ and $X:={\mathbb P}(\sF)$. As soon as $a_i>0$ for each index $i=2,\ldots, n-1$, the exact sequence \eqref{lucianex} expresses $Y$ as a smooth ample divisor of $X$; it is recovered by pushing forward the exact sequence
\[0\to \sO_X\to \sO_X(Y)\to \sO_Y(Y)\to 0.\]
\end{rem*}

The following fact is well known. We include the proof for reader's convenience.

\begin{lemma}\label{wk} Let $V={\mathbb P}(\sE)$ be a $\pn {n-1}$-bundle over $\pn 1$, for some  rank $n$ vector bundle $\sE=\bigoplus_{i=1}^n\sO_{\pn 1}(a_i)   $. Assume that $V$ is a Fano manifold, of index $i(V)$. Then, for some integer $a$, either
\begin{enumerate}
\item[{\rm (i)}] $V= \pn {n-1}\times \pn 1$, $i(V)\leq 2$  and $\sE=  \bigoplus_{i=1}^n\sO_{\pn 1}(a)$; or
\item[{\rm (ii)}] $V$ is the blowing-up, $\grs:V\to \pn n$, along a codimension  two  linear subspace of $\pn n$, $i(V)=1$ and $\sE= \bigoplus_{i=1}^{n-1}\sO_{\pn 1}(a) \oplus\sO_{\pn 1}(a+1)$.\end{enumerate}
\end{lemma}
\proof
After normalization of the integers $a_i$, write
$0=a_1\leq a_2\leq \cdots\leq a_n$ and consider the section $\Gamma$ of $p:V\to \pn 1$ corresponding to the quotient 
\[\sE\to \sO_{\pn 1}\to 0\] onto the first summand. Note that the morphism $\vphi_{|\xi_V|}:V\to \pn m$ maps the curve 
$\Gamma$ to a point. On the other hand, since $V$ is a Fano manifold  with $\Pic V\cong \zed\oplus\zed$,  there are two extremal rays, $R_1$, corresponding to $p$, and $ R_2$, generating the cone $\overline{NE}(V)$. Since the morphism 
$\vphi_{|\xi_V|}$ is not finite, it must  coincide with the contraction  of $R_2$. Now, setting $d:=\sum_{i=1}^na_i$, the canonical bundle formula yields
\[-K_V\cong p^*\sO_{\pn 1}(2-d)\otimes \sO_V(n).\] 
Therefore, dotting with $\Gamma$, we get $d<2$, so that either $\sE\cong  
\bigoplus_{i=1}^n\sO_{\pn 1}$, or  $\sE\cong\sO_{\pn 1}\oplus \cdots \oplus \sO_{\pn 1}\oplus\ \sO_{\pn 1}(1)$, leading to the two cases as in the statement. Note that by the canonical bundle formula, in the first case   the index of $V$ is 
$i(V)\leq 2$, while, in the second case, $i(V)=1$.
\qed

\begin{examples*}\label{only}(Only known examples of non-extendable extremal rays). Let $X$ be a projective manifold of dimension $n\geq 4$. Let $Y$ be a smooth ample divisor on $X$. The only known examples of extremal rays $R$ of $Y$ which do not extend to $X$ are the following:
\begin{enumerate}
\item[(1)] $Y=\pn 1\times\pn {n-2}$, and $R$ is the nef extremal ray  corresponding to the $\pn 1$-bundle projection $q:Y\to \pn {n-2}$. The manifold  $X$ is constructed as in Remark~\ref{lurem};
\item[(2)] $Y={\mathbb P}\big(\bigoplus_{i=1}^{n-2}\sO_{\pn 1}(a) \oplus\sO_{\pn 1}(a+1)\big)$ for some integer $a\geq 2$,  $R$ is the non-nef extremal ray corresponding to the blowing-up, $\grs:Y\to \pn {n-1}$, along a codimension  two  linear subspace of $\pn {n-1}$. Again, $X$ is constructed as in Remark~\ref{lurem}, $X\neq \pn 1\times\pn {n-1}$.
\end{enumerate}
In case (1), the $\pn {n-2}$-bundle  projection $p:Y\to \pn 1$ on the first factor extends by construction. Then, if $q:Y\to \pn {n-2}$ extends too, we would have a surjective map $\pn {n-1}\to \pn {n-2}$, where $\pn {n-1}$ is a fiber of the extension of $p$; a contradiction.

Let  $Y$ be as in case (2). By Proposition~\ref{Badprop}, $Y$ embeds  as a smooth  ample divisor of $X:={\mathbb P}\big(\sO_{\pn 1}(1)\oplus\sO_{\pn 1}(a-1)\oplus
\sO_{\pn 1}(a)^{\oplus (n-3)} \oplus\sO_{\pn 1}(a+1)\big)$.  Note that  $\overline{NE}(Y)=\langle R_1,R_2\rangle$, where $R_1$, $R_2$ are   the extremal rays corresponding to the bundle projection $Y\to \pn 1$, and to the blowing-up $\grs:Y\to \pn {n-1}$ respectively.  Moreover, Lemma~\ref{wk} applies to say that $X$ is not a Fano manifold.  Therefore $\overline{NE}(Y)$ is strictly contained in $\overline{NE}(X)$.
Since $\varrho(Y)=\varrho(X)$ by the Lefschetz theorem,  and  the projection $Y\to \pn 1$ extends
by construction,  we thus conclude that the extremal ray $R_2$ does not extend to $X$.
Clearly $X\neq \pn 1\times\pn {n-1}$ in the above example. Note that by taking as $Y$ a hyperplane section of the Segre embeding $X$ of $\pn 1\times\pn {n-1}$,  the restriction to $Y$ of the bundle projection $X\to \pn {n-1}$ is in fact the blowing-up $\grs:Y\to \pn {n-1}$ along a codimension  two  linear subspace. Of course, in this case, the extremal ray defining $\grs$ extends to $X$. In terms of Proposition~\ref{Badprop}, this situation corresponds to the case when  $a=1$, that is $Y={\mathbb P}\big(\bigoplus_{i=1}^{n-2}\sO_{\pn 1}(1) \oplus\sO_{\pn 1}(2)\big)$ is an ample divisor of $X={\mathbb P}\big(\bigoplus_{i=1}^n\sO_{\pn 1}(1)\big)$.
\end{examples*}

That those in \ref{only} are the only known examples of non-extendable extremal rays looks quite surprising to us. We propose the following speculations with the hope they may eventually lead to an explanation of this fact.

Let  $X$ be a projective manifold of dimension $n\geq 4$. Let $Y$ be a smooth ample divisor on $X$.
 Consider   the following property,  for an extremal ray, $R$, of $Y$.
\begin{center} \noindent $(\star)\;\;\;\;$ {\em For any nef divisor $H$ on $Y$, such that 
$H^{\perp}\cap(\overline{NE}(Y)\setminus \{0\})=R$, it follows\\\hskip30pt that $\overline{H}$ is nef  {\rm (}here  $\overline{H}$ is the unique divisor class on $X$ such that $\overline{H}_{Y}=H)$}.
\end{center}

First, note that $(\star)$ implies that $R$ is  a ray of $X$ (cf. \cite[Theorem 4.1]{BdFL} and end of proof of Theorem~\ref{Occhetta}). 

\begin{prop}\label{sexy} Assume that property $(\star)$ holds for all extremal rays of $Y$ such that $R\subset \scd^{\perp}$ for some nef divisor $\scd$ on $X$, $\scd \not \sim 0$.
Assume, also, that $\varrho(X)\geq 3$. Then every extremal ray of $Y$ extends to an extremal ray of $X$.
\end{prop}
\proof  Assume that we have some extremal ray of $Y$, say $R_0$, which is  not a ray of $X$. We may assume that $R_0\subset \overline{NE}(X)_{(K_X+(1+\varepsilon)Y) <0}$ for some $\varepsilon >0$.

Observe that, by our hypothesis, any ray $R$ of $Y$ which is contained in some face of $\overline{NE}(X)$  satisfies property ($\star$), and hence, as noted above, $R$ is a ray of $X$.

Let $\scc_\circ:=(K_X+(1+\varepsilon)Y)^{\perp}\cap (\overline{NE}(Y)\setminus \{0\})$ and let $R_0,R_1,\ldots,R_s$ be all the extremal rays of $\overline{NE}(Y)$ such that $R_i\subset \overline{NE}(X)_{(K_X+(1+\varepsilon)Y)<0}$, $i=0,1,\ldots,s$. Next consider the non-degenerate convex cone $\scc\subset \reals^{\varrho(Y)}$ defined by
\[\scc=\langle \scc_{\circ},R_0,R_1,\ldots,R_s\rangle=\overline{NE}(Y)\cap \overline{NE}(X)_{(K_X+(1+\varepsilon)Y)\le 0}.\]
Then $\scc_{\circ}$ is a face of the cone $\scc$ and, for each index $i$,
 
 \[R_i\not\subset \langle R_0,R_1,\ldots,R_{i-1},R_{i+1},\ldots,R_s\rangle.\]

\begin{clam}\label{percepito} There exists a face $F$ of the cone $\scc$ such that $F$ contains two extremal rays of $Y$, say $R', R''$, such that $R'$ is a ray of $X$ and $R''$ is not. 
\end{clam}

\proof  Recall that by Proposition \ref{KollarProp} all extremal rays of $X$ contained in $\overline{NE}(X)_{(K_X+Y)\leq 0}$  are also extremal rays of $Y$. Therefore, some of the extremal rays $R_j$, $j\geq 1$, are rays of $X$. Take one of
them, say $R_1$. Then there exists a face $F_1$ of $\scc$ containing $R_1$ and some of the other
 rays $R_{j'}$,  for some $j'\in \{2,\ldots,s\}$. If one of the rays $R_{j'}$ (say $R_2$), is not a ray of $X$, then take $R'=R_1$, $R''=R_2$ and $F=F_1$.
 If all the extremal rays $R_{j'}$,  $j'\in \{2,\ldots,s\}$, are rays of $X$, we apply the same argument to conclude that either we prove the claim, or every extremal ray of $Y$ contained in $\scc$ lifts to an extremal ray of $X$, contradicting the assumption that $R_0$ does not. 
 
 Thus we may assume to be in the situation described in Claim \ref{percepito}.
 Take a nef divisor $D$ on $Y$ such that $F=D^{\perp}\cap(\overline{NE}(Y)\setminus\{0\})$. If the unique divisor class $\overline{D}$ on $X$ which restricts to $D$ is nef, then by our assumption $R''$ is a ray of $\overline{NE}(X)$. This  contradicts the claim, so $\overline{D}$ is not nef on $X$.

 Now, take a nef divisor $\sch$ on $X$ such that $R'=\sch^{\perp}\cap(\overline{NE}(X)\setminus\{0\})$. For $0\leq \gra\leq 1$, consider $D_\gra=\gra\overline{D}+(1-\gra)\sch$ and let $\lambda:=\sup\{\gra\;|\; D_\gra\; {\rm is \; nef}\}$. Then $0\leq \lambda<1$.
 
\vskip-1.5cm
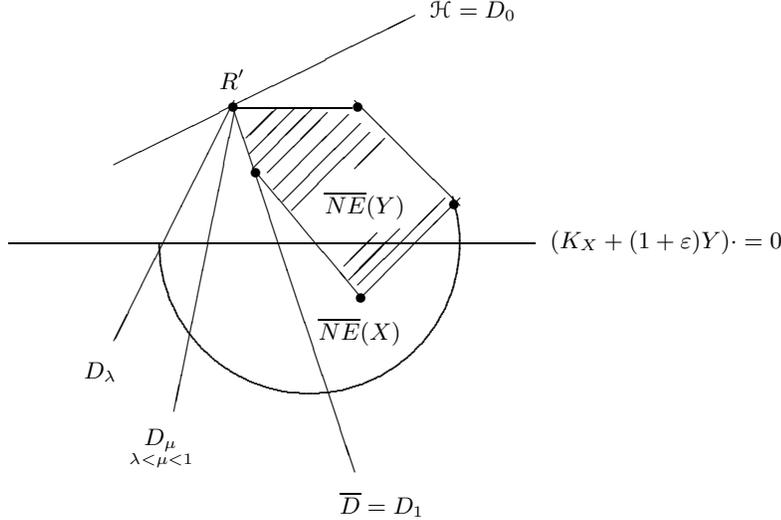
\begin{figure}[h]\label{nuco} 
{\fontsize{9pt}{12pt}
\setlength{\unitlength}{2cm} 
\hskip2cm\begin{center}\begin{picture}(2,2) 
\arc(-1,0){198}
\put(-0.42,0,7){$\line(1,1){0.2}$}
\put(-0.4,0,6){$\line(1,1){0.3}$}
\put(-0.38,0,5){$\line(1,1){0.4}$}
\put(-0.28,0,47){$\line(1,1){0.4}$}
\put(0.1,0.2){$\hbox{$\overline{NE}(Y)$}$}
\put(0.06,-0.64){$\hbox{$\overline{NE}(X)$}$}
\put(-0.24,0.36){$\line(1,1){0.5}$}
\put(-0.18,0.28){$\line(1,1){0.52}$}
\put(-0.12,0.22){$\line(1,1){0.52}$}
\put(0.3,-0.2){$\line(1,1){0.3}$}
\put(0.32,-0.28){$\line(1,1){0.58}$}
\put(0.26,-0.24){$\line(1,1){0.25}$}
\put(0.23,-0.15){$\line(1,1){0.22}$}
\put(0.3,0.5){$\line(1,1){0.2}$}

\put(0.92,0.22){$\bullet$}
\put(1,0.25){$\line(-1,1){0.7}$}
\put(0.28,0.87){$\bullet$}
\put(0.28,0.9){$\line(-1,0){0.8}$}
\put(-0.55,0.87){$\bullet$}
\put(-0.5,0.88){$\line(1,-3){0.8}$}
\put(-0.5,0.95){$\line(-1,-2){0.8}$}
\put(-0.4,0.43){$\bullet$}
\put(-0.35,0.47){$\line(5,-6){0.7}$}
\put(-2,0){$\line(1,0){3.5}$}
\put(1.6,-0.03){$(K_X + (1+\varepsilon)Y)\cdot = 0$}
\put(-0.6,1.03){$R'$}
\put(0.8,1.5){$\sch=D_0$}
\put(-1.5,-0.9){$D_{\lambda}$}
\put(-1.2,-1.4){$\substack{{\displaystyle D_{\mu}} \\ \lambda < \mu <1}$}
\put(0.2,-1.8){$\overline{D} = D_1$}
\put(-0.5,0.88){$\line(-1,-5){0.4}$}
\put(-1.3,0.52){$\line(2,1){2}$}
\put(0.3,-0.4){$\bullet$}
\put(0.32,-0.38){$\line(1,1){0.68}$}
\end{picture} \end{center}}
\vskip3.4cm 
\caption{The two cones of curves. }\end{figure}
\begin{clam}\label{notprop} The extremal ray $R'$ does not satisfy property $(\star)$. 
\end{clam} 
\proof For some $\lambda<\mu<1$, the divisor $D_\mu$ is not nef and its restriction  ${D_\mu}_{|Y}=\mu D+(1-\mu)\sch_Y$ is nef on $Y$. In particular one has
\[R'\subseteq ({{D_\mu}_{|Y}})^{\perp}\cap(\overline{NE}(Y)\setminus\{0\}) \subseteq \sch_Y^{\perp}\cap(\overline{NE}(Y)\setminus\{0\}).\]
Thus, since $R'=\sch^{\perp}\cap(\overline{NE}(X)\setminus\{0\})$,  equalities hold above, so that  $R' =({{D_\mu}_{|Y}})^{\perp}\cap(\overline{NE}(Y)\setminus\{0\})$. Thus 
  $R'$ does not satisfy property $(\star)$; see Figure 1.

 The  claim above leads to a contradiction, and hence concludes the proof of the proposition.
\qed
\begin{rem*} (1) The examples of non-extending extremal rays in Example~\ref{only} were known to the experts since the early 80's. The 
resulting observation that, in general, $\overline{NE}(Y)\subsetneq \overline{NE}(X)$ when $Y\subset X$ is ample, was rediscovered
in \cite{HLW}.

(2) The (open) problem of deciding whether or not the property ($\star$) holds for all extremal rays $R$ of $Y$ such that 
$R\subset D^\perp$ for  some nef divisor $D$ on $X$, $D\not \sim 0$, looks hard but very interesting. A positive
answer to it would imply, via Proposition~\ref{sexy}, that non-extendable rays may only occur when $\varrho(X) =2$. Moreover, a positive 
answer would give a proof for the claim made in \cite[Theorem~4.3]{HLW}. See also Remark~\ref{beautyfulmind} for further applications. 
\end{rem*}

\section{Applications to $\pn d$-bundles  and blowing-ups}\label{application}\addtocounter{subsection}{1}\setcounter{theorem}{0}

Let us start  by  recalling some useful preliminary facts.

  The following result goes back to Goren \cite{Gor} and  Kobayashi--Ochiai \cite{KoOc}.  We refer also to Fujita  \cite[Chapter I, (1.1), (1.2)]{FuBook} where the Cohen--Macaulay assumption was removed.

\begin{theorem}\label{KOF} Let $L$ be an ample line bundle on an $n$-dimensional irreducible projective variety $X$. If $L^n=1$ and $h^0(L)\geq n+1$, then $(X,L)\cong (\pn n, \sO_{\pn n}(1))$.
\end{theorem}

\begin{corollary}\label{Co} If $X$ is a Fano manifold of dimension $n$ and index $i$, we have $i\le n+1$ and $X\cong \pn n$ if equality holds.
\end{corollary}\proof Use the Hilbert polynomial and Kodaira vanishing to check the hypothesis of Theorem~\ref{KOF}. See \cite{KoOc} for details.
\qed

\begin{rem*}\label{Re} In a completely similar way one proves that if $X$ is as above and $i=n$, then $X\cong \sQ^n$, see again \cite{KoOc}.
\end{rem*}
The smooth version of  the following fact was proved  by Ramanujam \cite{Ramanan} and Sommese \cite{So1}.  The general case is due to B\v adescu \cite{BadHD}. We also refer to   \cite[Section~2.6 and (5.4.10)]{Book} for a different argument based on Rossi's extension theorem.

\begin{theorem}\label{Pn} Let $Y\cong \pn {n-1}$ be an ample  Cartier divisor  on a normal projective variety  $X$ of dimension $n\geq 3$.  Then $X$ is the cone $\sC(\pn n,\sO_{\pn n}(s))$, where  $\sO _Y(Y)\cong \sO_{\pn {n-1}}(s)$. If $s=1$, the assertion is true for $n=2$ as well. \end{theorem}\proof Let us recall the argument in the smooth case.

First, note that,  by Lefschetz theorem, $i^*:\Pic X\stackrel{\sim}{\to}\Pic Y$ via the embedding $i: Y\hookrightarrow  X$. This is clear if $n\geq 4$. If $n=3$, then $i^*:\Pic X\to \Pic Y$  is injective with torsion free cokernel. Therefore $i^*$ is an isomorphism, since $\Pic Y\cong \zed$.

Let $L$ be the ample generator of $\Pic X$ and assume $Y \in |aL|$ for some $a\geq 1$. Note that $L_Y = \sO _{\pn{n-1}} (1)$. It then follows that
\[1=(L_Y)^{n-1}_Y=(L^{n-1}\cdot Y)_X= a(L^n).\]
Therefore $a=1$ and $L^n=1$.

Next, Kodaira vanishing and the exact sequence
\[0\to \sO_X(-L)\to \sO_X\to \sO_Y\to 0\]
 give $H^1(X,\sO_X)=0$. Hence the exact sequence
\[0\to \sO_X\to L\to L_Y\to 0\] yields
$h^0(L)= h^0(L_Y)+1=\dim X+1$. Then Theorem~\ref{KOF} applies to give $(X,L)\cong
 (\pn n, \sO_{\pn n}(1))$.
\qed

For the first statement of the theorem below, see Fujita \cite[(2.12)]{FuSemi} or Ionescu  \cite[p.\ 467]{IoGen}; the  second point follows from  Theorem~\ref{Pn} and some of   B\u adescu's arguments  in \cite{BadAmp1, BadAmp2}; the third point was noticed in 
\cite[Section 2]{BSCe}  and \cite{BSW}.

\begin{theorem}\label{FB} Let $X$ be an $n$-dimensional projective manifold.
\begin{enumerate}
\item[{\rm (i)}] Let $\pi:X\to Z$ be a surjective morphism from $X$ onto a normal variety $Z$. Let $L$ be an ample line bundle on $X$. Assume that $(F,L_F)\cong (\pn d,\sO_{\pn d}(1))$ for a general fiber, $F$, of $\pi$ and that all fibers of $\pi$ are $d$-dimensional. Then $\pi:X\to Z$ is a linear $\pn d$-bundle with $X={\mathbb P}(\pi_*L)$.
\item[{\rm (ii)}] Let $p:Y\to Z$ be a $\pn d$-bundle over a projective manifold  $Z$. Assume that $Y$ is an ample divisor on $X$. Furthermore assume that $p$ extends to a morphism  $\overline{p}:X\to Z$. Then there exists a non-splitting exact sequence
\[0\to \sO_Z\to \overline{\sE}\stackrel{u}{\to} \sE\to 0,\]
where $\overline{\sE}$, $\sE$ are ample vector bundles on $Z$
such that $X={\mathbb P}(\overline{\sE})$, $Y={\mathbb P}(\sE)$, $p$,  $\overline{p}$ are the bundle projections on $Z$, and the inclusion $Y\subset X$ is induced by $u$.
\item[{\rm (iii)}] Let $\pi:X\to Z$ be a linear $\pn {d+1}$-bundle over a projective manifold $Z$. Assume that $\dim Z<d+1$ and the tautological line bundle of $X$, say $L$, is ample. Then $K_X+(d+2)L$ is nef. 
Moreover, the bundle projection $\pi$ is associated to the linear system $|m(K_X+(d+2)L)|$ for $m\gg 0$ {\rm (}i.e., $(X,L)$ is a scroll over $Z${\rm )}.
\end{enumerate}
\end{theorem}
\proof (i)  Following  the argument  as in  \cite[p.\ 467]{IoGen}, let us first show that $Z$ is smooth. Indeed, let $z\in Z$ be a closed point and denote by $\Delta$  the fiber over $z$. Consider the embedding of $X$ given by $|mL|$ for $m\gg 0$. Let $\widetilde{Z}$ be the smooth $(n-d)$-fold got by intersecting $d$ general members $H_1,\ldots,H_d$ of $|mL|$. Since $(F,L_F)\cong (\pn d,\sO_{\pn d}(1))$ for a general fiber  $F$,  the restriction $\widetilde{p}$ of $p$ to $\widetilde{Z}$ has degree $m^d$.
Since $p$ has equidimensional fibers,  $\widetilde{Z}\cap \Delta$ is a $0$-dimensional scheme, and its length   is given by $\ell(\widetilde{Z}\cap \Delta)=m^d(L_\Delta^d)\geq m^d$. Since $Z$ is normal, it follows by a well known criterion (see e.g., \cite[Chapter II, Theorem~6]{Sh}) that  the above inequality is in fact an equality. Hence in particular $L_\Delta^d=1$, so that $\Delta$ is irreducible and generically reduced.  Therefore, by the generality of $H_1,\ldots,H_d$,  we may assume that  
$\widetilde{Z}\cap \Delta$ is  a reduced $0$-cycle  consisting of $\#(\widetilde{Z}\cap \Delta)=\ell(\widetilde{Z}\cap \Delta)=\deg (\widetilde{p})$ distinct points. It thus follows that $\widetilde{p}$ is \'etale over $z$. Therefore $Z$ is smooth at $z$ since $\widetilde{Z}$ is smooth.

 Since all fibers of 
$\pi$  are equidimensional and $X$ and $Z$ are  smooth, the morphism $\pi$ is flat. Let now $\Delta$ be any fiber of $\pi$. We have seen above that $\Delta$ is irreducible and generically 
reduced.
Moreover, $\Delta$ is Cohen--Macaulay since every fiber is defined by exactly $n-d$
coordinate functions. It thus follows that $\Delta$ is in fact reduced. By the
semicontinuity theorem for  dimensions of spaces of sections on fibers of a flat
morphism \cite[Chapter III, Theorem~12.8]{Hrt},  $h^0(L_\Delta) \ge h^0(L_F)=d+1$. Then by   Theorem \ref{KOF}  we conclude that
$(\Delta,L_\Delta) \cong  (\pn {d},
\sO_{\pn {d}}(1))$ for every fiber $\Delta$ of $\pi$. Then $X\cong{\mathbb P}(\sE)$, where $\sE := \pi_*L$.

(ii) Note that, by ampleness of $Y$, $p$ equidimensional implies $\overline{p}$ equidimensional.  Let $F$ be a general fiber of $\overline{p}$ and let $f=F\cap Y$ be the corresponding fiber of $p$. By (i), it is enough to show that $F\cong \pn {d+1}$ and $L_F\cong\sO_{\pn {d+1}}(1)$, where $L=\sO_X(Y)$. If $d \geq 2$ we conclude by Theorem \ref{Pn}.

Thus we can assume $d=1$. By taking general hyperplane sections of $Z$ and by base change, we can also assume  that $Z$ is a smooth curve (and hence $F$ is a divisor). Here we follow B\u adescu's argument. By the Lefschetz theorem, we have an exact sequence
\[0\to {\rm Num}(X)\stackrel{i^*}{\to}{\rm Num}(Y)\to {\rm Coker}(i^*)\to 0,\]
 where $i^*$ is the morphism induced by the embedding $i:Y\hookrightarrow X$ and ${\rm Coker}(i^*)$ is torsion free. 
First note that   ${\rm Num}(X)\not\cong \zed$, since otherwise $F$ would be an ample divisor. Since 
${\rm Num}(Y)\cong\zed\oplus\zed$, we thus conclude that  ${\rm Num}(X)\cong\zed\oplus\zed$, and therefore that ${\rm Coker}(i^*)=(0)$  since it is torsion free. Let $h$ be a section of the bundle $p:Y\to Z$. Then ${\rm Num}(Y)=\zed \langle f\rangle\oplus\zed \langle h\rangle$ and 
${\rm Num}(X)=\zed \langle F\rangle \oplus\zed \langle H\rangle $ for some line bundle $H$ on $X$ inducing $h$ on $Y$. Write $Y\sim aF+bH$ for some integers $a$, $b$. Since $h\cdot f=H\cdot Y\cdot F=1$, we get $1=b(H^2\cdot F)$ and therefore $b=\pm 1$. Then $Y\cdot f=(aF\pm H)\cdot f=\pm(H\cdot f)=\pm 1$. By ampleness of $Y$, it must be $b=1$. Thus $f\cong \pn 1$ has self-intersection $f^2=Y^2\cdot F=Y\cdot f=1$ on $F$. This implies $F\cong\pn 2$  and $L_F\cong\sO_{\pn 2}(1)$ by using Theorem \ref{Pn}. Applying the first part of the statement and   pushing forward under $\overline{p}$ the exact sequence 
\[0\to \sO_X\to \sO_X(Y)\to \sO_Y(Y)\to 0\] we get the desired conclusion.

(iii)  If $K_X +(d+2)L$ is not nef, there exists an extremal ray $R$ of $X$ such that $(K_X+(d+2)L)\cdot R<0$. It follows that $\ell(R)>d+2$. Therefore,  if $\Delta$ is a positive dimensional fiber of the contraction ${\rm cont}_R$, we have $\dim\Delta\geq\ell(R)-1\geq d+2$ (see Theorem~\ref{ExRayIW}).
Then, for a fiber $F$ of $\pi$, we have
\[\dim F+\dim \Delta\geq 2d+3 >n+1\]
(where the last inequality follows from the assumption $d+1>\dim Z$, which is equivalent to saying that $2d+2> n$). Hence $\dim(F\cap \Delta)\geq 2$. Thus there exists a curve $C\subset F$ such that $(K_X+(d+2)L)\cdot C<0$, contradicting $(K_X+(d+2)L)_{|F}\approx 0$. Therefore we conclude that $K_X+(d+2)L$ is nef, and hence, by Theorem \ref{KRS}, the linear system $|m(K_X+(d+2)L)|$ defines a morphism, say $\vphi$, for $m\gg 0$.

Let now $R\subset (K_X+(d+2)L)^{\perp}\cap (\overline{NE}(X)\setminus\{0\})$ be an extremal ray. We have $(K_X+(d+2)L)\cdot R=0$, so that $\ell(R)\geq d+2$, and hence, as above,
$\dim\Delta\geq \ell(R)-1\geq d+1$.

Since $2d+2>n$, it follows, again by Theorem~\ref{ExRayIW}, that
 $\dim F+\dim \Delta\geq 2d+2\geq n+1$. Then $\dim(F\cap \Delta)\geq1$. This implies that $\pi$ is the contraction $ {\rm cont}_R$ of the extremal ray $R$. Since this is true for each extremal ray as above, we conclude that the face
 $(K_X+(d+2)L)^{\perp}\cap (\overline{NE}(X)\setminus \{0\})$  is in fact $1$-dimensional and that $p$ coincides with the  morphism $\vphi$.
\qed 

\begin{rem*}\label{inpart} (1) In the  boundary case $d+1=\dim Z$ of Theorem~\ref{FB}(iii), the same argument gives the nefness of  $K_X+(d+2)L$; 
moreover, further considerations show that  the bundle projection $\pi$ is associated to $|m(K_X+(d+2)L)|$ for $m\gg 0$   unless $X\cong\pn {d+1}\times\pn {d+1}$. We refer for this to 
\cite[(3.1)]{BSW}.

(2) Note that by Theorem~\ref{basere} or Lemma~\ref{easy} below,  statement  (iii) of Theorem~\ref{FB} applies under the conditions  in  \ref{FB}(ii).
\end{rem*}

In the case of  $\pn d$-bundles, Sommese's theorem~\ref{basere}  admits the following simple 
alternative proof.

\begin{lemma}\label{easy} Let $Y$ be a smooth  ample  divisor on
a projective manifold, $X$. Assume that $Y$ is a $\pn d$-bundle over a manifold $Z$. Further assume that $p$ has an extension $\overline{p}:X\to Z$. Then $\dim Z\leq d+1$.
\end{lemma}
\proof By Theorem~\ref{FB}
we get the exact sequence
\begin{equation}\label{E}
0\to \sO_Z\to \overline{\sE}\to \sE\to 0,
\end{equation}
where $\overline{\sE}$, $\sE$ are ample vector bundles  on $Z$. Arguing by contradiction,  assume $\dim Z>{\rm rk}\sE=d+1$, that is $1\leq \dim Z-{\rm rk}\sE$. By le Potier's vanishing theorem \cite{LeP} we have
$H^i(Z,\sE^*)=0$ for $i\leq \dim Z-{\rm rk}\sE$. Therefore $H^1(Z,\sE^*)=0$, so we conclude that (\ref{E})
splits, contradicting ampleness of $\overline{\sE}$.
\qed
\begin{lemma}\label{key1} Let $L$ be an  ample line bundle on a
 projective manifold, $X$, of dimension $n\ge 4$. Assume that there is a
smooth $Y\in |L|$ such that $Y$ is a $\pn d$-bundle, $p:Y\to Z$, over a  manifold $Z$.
Let $\ell$ be a line in a fiber $\pn d$ of $p$. Further assume that $H\cdot \ell=1$ for some ample line bundle $H$ on $X$.
Then  $(X,L)\cong ({\mathbb P}(\sE),\xi_{\mathbb P})$ for an ample rank $d+2$ vector bundle, $\sE$, 
on $Z$ with $p$ equal to the
restriction to $Y$ of the induced projection ${\mathbb P}(\sE)\to Z$.\end{lemma}

\proof 
Lines in the fibers of $p$ define a covering family of $Y$. By our assumptions, the induced family on $X$ is unsplit. 
Therefore Theorem~\ref{Occhetta} applies to give that $p$ extends. We conclude by Theorem~\ref{FB}(ii).
\qed

The following gives a precise answer to  Question \ref{zeroQ}.

\begin{prop}\label{Qnot}
Let $X$ be an $n$-dimensional projective manifold embedded in~$\pn N$, $n\geq 4$.
\begin{itemize}
\item[{\rm (i)}] Assume that $X$ has a smooth hyperplane section $Y= X\cap H$ which is a $\pn d$-bundle over a manifold $Z$,
say $p: Y \to Z$ , such that the fibers of $p$ are linear subspaces of $\pn N$. Then $p$ lifts to a linear $\pn{d+1}$-bundle
$\overline p : X\to Z$. Moreover, this is possible only if $d+1\geq \dim Z$.

\item[{\rm (ii)}] Conversely, assume that $\pi: X\to Z$ is a $\pn{d+1}$-bundle with linear fibers and 
$d+1 \geq \dim Z$. Then there exists a smooth hyperplane section $Y= X\cap H$ which is a $\pn d$-bundle.

\end{itemize}
\end{prop}
\proof (i) If $\ell$ is a line contained in some fiber of $p$, we have $H\cdot \ell=1$. Thus the first assertion follows from Lemma~\ref{key1} and the second from Theorem~\ref{basere} (or Lemma~\ref{easy}).

(ii) Consider the incidence relation
\[W:= \{ (z,h) \mid H\supseteq F_z\} \subseteq Z\times (\pn N)^{\vee},\]
where $F_z= \pi^{-1}(z)$ is the fiber $\pn {d+1}$ over a point $z\in Z$ and $h\in (\pn N)^{\vee}$ is the point corresponding to the hyperplane $H$ in $\pn N$. Then $\dim W = \dim Z +N -(d+1) -1$, so that $\dim Z \le d+1$ gives $\dim W \le N-1$. Therefore there exists a hyperplane $H$ in $\pn N$ not containing fibers of $\pi$ and we are done.  
\qed

Consider the setting as in diagram (\ref{diagram}) from Section~2 with $\dim X\geq 4$.  Let us discuss  here some  applications under the assumption that  the canonical bundle $K_Z$ of $Z$ is nef. We follow the exposition  in \cite[Section 4]{IoNef}, where the results are proved for a strictly nef and big divisor $Y$ on $X$.

As a first application, we consider the case when the morphism $p:Y\to Z$ is a $\pn d$-bundle.  
An  analogous result, assuming  $\grk(Z)\geq 0$ instead of  $K_Z$ to be nef, was proved in \cite{FS} in a completely different way.

The following result is essentially due to Wi\'sniewski \cite[(3.3)]{Wcontr}.

\begin{lemma}\label{PalUnique} Let $Y$ be a $\pn d$-bundle over a smooth projective variety and let $p:Y\to Z$ be the bundle
projection. If $K_Z$ is nef, then $Y$ admits a unique extremal ray,  and  $p$
is its contraction.
\end{lemma}
\proof  Assume by contradiction that there exists an extremal
rational curve, $C$, which is not contracted by
$p$. Let
$\Gamma\cong \pn 1$ be the normalization of $p(C)$. Let $f:\Gamma\to Z$ be the induced morphism and consider the base change diagram
\[ \xymatrix{
Y' \ar[r]^(.50)g \ar[d]_{p'} & Y \ar[d]^p \\
\Gamma \ar[r]_{f} & Z.} \]
Since $\Gamma$ is a smooth curve, we have a vector bundle $\sE$ on $\Gamma$, of rank $r:=d+1$, such that 
$Y'={\mathbb P}(\sE)$.
Let $F'$ be the fiber of the bundle projection $p'$ and let $C'\subset Y'$ be a curve mapped onto $C$ under $g$.
Let $\sT'$ be the tautological line bundle  on $Y'$. Then we get
\[0>(C\cdot K_Y)=(C'\cdot g^*K_Y) =-r(C' \cdot \sT')+ (C'\cdot(f\circ p')^*K_Z)+(C' \cdot
{p'}^*\det(\sE)).  \]
Since $K_Z$ is nef, it thus follows $r(C'\cdot \sT')>(C'\cdot{p'}^*\det(\sE))$. By the Grothendieck
theorem, we have $\sE\cong\bigoplus_{i=1}^r\sO_{\pn 1}(a_i)$, where $a_1\geq\cdots \geq a_r$. Then the
inequality above yields $r(C'\cdot \sT')>ra_r(C'\cdot F')$. Thus
\begin{equation}\label{P1}
(C'\cdot \sT')>a_r(C'\cdot F').
\end{equation}
Let $R=\reals_+[C]$ be the extremal ray generated by $C$. The composition  $\vphi:={\rm cont}_R\circ g$ is a
morphism  defined by a linear sub-system of $|\gra \sT'+\grb F'|$, for some $\gra>0$. Since  $C'\cdot (\gra
\sT'+\grb F')=0$, we get $\gra(C'\cdot \sT')=-\grb(C'\cdot F')$. Therefore (\ref{P1})  gives
\begin{equation}\label{P2}
-\grb>\gra a_r.
\end{equation}
Let $C_r$ be the section of $p'$ corresponding to the surjection  $\sE \to \sO_{\pn 1}(a_r)\to 0$. Then
$C_r\cdot (\gra\sT'+\grb F')\geq 0$, which contradicts (\ref{P2}).
\qed

\begin{prop}\label{Lucian1} Let $p:Y\to Z$ be a $\pn d$-bundle over a projective manifold $Z$. Assume that  $K_Z$ is nef.   If $Y$ is an ample divisor on a manifold $X$, then $p$ extends to a morphism $\overline{p}:X\to Z$.  Furthermore  there exists a non-splitting  exact sequence
\[0\to \sO_Z\to\overline{\sE}\to\sE\to 0\] such that $X={\mathbb P}(\overline{\sE})$, $Y={\mathbb P}({\sE})$, and $p$, $\overline{p}$ are the bundle projections on $Z$.
\end{prop}
\proof By Lemma \ref{PalUnique}, the bundle projection $p$ is the contraction of the unique extremal ray on $Y$. Then, by Corollary~\ref{newcor}, $p$  extends.  Thus Theorem~\ref{FB}(ii) applies to give the result.
\qed

Next, we consider the case when the morphism $p:Y\to Z$ is a blowing-up. The following general fact is a direct consequence of  Lemma~\ref{PalUnique}.

\begin{lemma}\label{Palbup1} Let $Z$ be a projective manifold. Let $p:Y\to Z$ be the blowing-up along a smooth subvariety $T$ of codimension $\geq 2$. Assume that the canonical bundles $K_T$ and $K_Z$ are nef.  Then $Y$ has only one extremal ray.
\end{lemma}
\proof  Let $c$ be the codimension of $T$ in $Z$, and let $E$ be the exceptional divisor of $p$. Let $R=\reals_+[C]$ be  any extremal ray on $Y$. Since $K_Y\approx p^*K_Z\otimes \sO_Y((c-1)E)$, we conclude that $E\cdot C<0$. Therefore  $C$ is contained in $E$ and $K_E\cdot C<0$. Then apply the proof of Lemma \ref{PalUnique} to the $\pn {c-1}$-bundle $E\to T$.\qed

The following generalizes a result due to Sommese concerning the  reduction map in the case of threefolds (see \cite{Somin}, \cite[Theorem~I]{SoVar},  and  also \cite{IoDeg,IoGen}) and is closely related to Fujita's results in \cite{FuLef}.

\begin{prop}\label{Palbup2} Let $Z$ be a projective manifold. Assume that $K_Z$ is nef. Let $p:Y\to Z$ be the blowing-up along a smooth subvariety $T$ of codimension $c\geq 2$, such that $K_T$ is nef. If $Y$ is an ample divisor on a manifold $X$, then there exists a commutative diagram  \eqref{diagram}, where $W$ is  a smooth projective variety and either
\begin{enumerate}
\item[{\rm (i)}] $\overline{p}:X\to W$ is the blowing-up of $W$ along the image of $T$. Moreover,  $\overline{p}(Y)$ is an ample divisor on $W$ whenever $T$is $0$-dimensional; or
\item[{\rm (ii)}] $X$ is generically a $\pn 1$-bundle over $Z$  and $Y$ is a rational section of it. Moreover, $\dim \overline p^{-1}(T)=n-2$, $c=2$ and fibers of $\overline{p}$ are at most two-dimensional.\end{enumerate}
\end{prop} 
\proof By Lemma \ref{Palbup1}, $Y$ contains a unique extremal ray and $p$ is its contraction. Then, by 
Corollary~\ref{newcor}, $p$  extends to a contraction $\overline{p}:X\to W$ of an extremal ray on $X$, which gives rise to a commutative diagram (\ref{diagram}). Assume first that $\overline p$ is birational.
Let $E$, $\overline{E}$ be the exceptional loci of $p$, $\overline{p}$ respectively, so that $E=\overline{E}\cap Y$. As the restriction $p_{|E}:E\to T$ of $p$ to $E$ is  a $\pn {c-1}$-bundle, it follows from Theorem~\ref{FB}(ii) that $\overline{p}_{|\overline{E}}:\overline{E}\to T$ is a $\pn c$-bundle. 
It is now standard to see that $W$ is smooth, $Z$ is contained in $W$ as a divisor and $\overline p$ is the blowing-up of $W$ along $T$, cf. also \cite[Section 5]{FuLef}. In \cite[Section 5]{FuLef}  it is also proved that $\overline{p}(Y)\cong Z$ is an ample divisor on $W$ under the extra assumption that the restriction to $T$ of the line bundle $\sO_W(Z)$ is ample. In particular, $\overline{p}(Y)$ is ample on $W$ if $T$ is $0$-dimensional.

Now, assume that $\overline p$ is not birational. Then $\alpha:Z\to W$ from (\ref{diagram}) is an isomorphism and $Y$ is a rational section for $\overline p$, which is generically a $\pn 1$-bundle. 
Let $t\in T$ be a general point and let $l\subset F:=p^{-1}(t)\cong \pn {c-1}$ be a line. We put $a:= Y\cdot l$ and we denote by $f$ a general fiber of $\overline p$. Since $l$ is contracted by $\overline p$, it is numerically proportional to $f$. It follows easily that $l\sim af$ as 1-cycles. As $K_X\cdot f=-2$, we get $-2a=K_X\cdot l=K_Y\cdot l-Y\cdot l=1-c-a$. So $a=c-1$. Assume that $\dim \overline p^{-1}(T)=n-1$. We obtain that $0=a(\overline p^{-1}(T)\cdot f)=\overline p^{-1}(T)\cdot l =-1$, a contradiction. So $\dim \overline p^{-1}(T)=n-2$. Now, denote by $V$ the family of all deformations of $l$ in $X$.
We claim that $\dim \mbox{Locus}(V)\geq n-1$. Assuming the contrary, we would have $\dim \mbox{Locus}(V)\leq n-2$. From the exact sequence 
\[0\to N_{l/Y}\to N_{l/X}\to \mathcal{O}_l(a)\to 0,\]
using standard facts from deformation theory of rational curves, we find that $\dim V=h^0(N_{l/X})=n+c+a-4$. Thus, if $x$ is a point on $l$, \[\dim V_x\geq \dim V+1-(n-2)=c+a-1.\] But the same exact sequence gives \[\dim V_x\leq h^0(N_{l/X}(-1))\leq c+a-2.\] This is a contradiction and the claim is proved.
Since $\dim \overline p^{-1}(T)=n-2$, the claim implies that some deformation of $l$ equals a fiber of $\overline p$. In particular, $a=1$, so $c=2$ and we are done.
\qed

\begin{rem*}\label{Furem}  Let us explicitly point out that an analogous result  was proved by  Fujita \cite[Section 5]{FuLef}, under  the assumption that ${\rm codim}_ZT\geq 3$, but with no nefness condition on $K_T$ and $K_Z$. However, the contraction morphism $X\to W$ obtained in \cite{FuLef} is in general analytic, not necessarily projective. 
\end{rem*}

Next, let us consider the case when $Y$ admits a  pluricanonical  fibration. 

Recall that a {\em Calabi--Yau  manifold}\/ $Y$  is a projective variety  with
trivial canonical bundle and  $H^i(Y,\sO_Y)=0$ for $i=1,\ldots,\dim Y-1$.

\begin{prop}\label{IoCY} Let $X$ be a projective manifold of dimension $\geq 3$.  Let $Y$ be a
smooth ample divisor on $X$.  Assume that $K_Y$ is nef. Then the linear system $|m(K_X+Y)|$ is base points free for $m\gg 0$. If $K_Y$ is numerically trivial, then $Y$ is a Calabi--Yau variety and $X$ is a Fano manifold. If $(K_Y)^{k+1}$ is a trivial cycle and $(K_Y)^k$ is non-trivial for $0< k< n-1$, then $X$ is a Fano fibration over $Z$ and $k=\dim Z$. \end{prop}
\proof The proof runs parallel to that of Theorem ~\ref{six}. Since $K_Y$ is nef and $Y$ is ample,  we conclude that $K_X+Y$ is nef. Thus $m(K_X+Y)$ is spanned for $m\gg 0$ by the Kawamata--Reid--Shokurov base point  free theorem, and it defines  a morphism $\pi:X\to W$. By restricting to $Y$ we find that $|mK_Y|$ is base points free for $m\gg 0$. Then $Y$ admits a pluricanonical map, say $\vphi:=\vphi_{|mK_Y|}$.

If $(K_Y)^{n-1}= 0$, the morphism  $\vphi$ is a fibration. If $K_Y$ is numerically trivial, then  $K_X+Y$ is also, and thus $X$ is a Fano manifold. Hence in particular $H^i(X,\sO_X)=H^i(Y,\sO_Y)=0$ for $0<i<\dim Y$, so that $Y$ is a Calabi--Yau manifold. The remaining part of the statement is  clear.
\qed

We say that a projective manifold $Y$ is {\em extendable} if there exists a projective manifold $X$ such that $Y\subset X$ is an ample divisor. 
Proposition~\ref{IoCY} shows that a manifold $Y$ such that $K_Y$ is numerically trivial and either $K_Y$ is not  linearly trivial  or $h^1(\sO_Y)>0$, e.g., $Y$ an abelian variety, is not extendable.

It is worth noting that  Proposition~\ref{IoCY} shows that the Abundance conjecture \cite{KMM} holds true for extendable manifolds. Let us recall  what the conjecture says in the smooth case. {\em Let $Y$ be a projective manifold with $K_Y$ nef. Then $mK_Y$ is spanned by its global sections for $m\gg 0$}. 

We refer to \cite{IoNef} for a  further discussion  in the case when $Y$ is a  strictly nef divisor on $X$.

\section{Complete results in  the three dimensional case}\label{n3}\addtocounter{subsection}{1}\setcounter{theorem}{0}

Throughout this section we assume that $X$ is a smooth projective three fold and $Y \subset X$ is a smooth ample divisor. The 
following theorem implies a number of results from \cite{BadAmp1,BadAmp2, BadAmp3, Somin,SoVar} and \cite{IoDeg}. We follow the arguments in 
\cite{IoGen}. Note that \cite{Mo1} contains a precise description of all types of extremal rays of $X$.

\begin{definition}\label{6.1} (cf.\ \cite{Somin}) A {\it reduction} of the pair $(X,Y)$ is another pair $(X',Y')$, with $Y'\subset X'$ a smooth
ample divisor, such that $X'$ is got by  contracting all $(-1)$ planes $E \cong \mathbb{P}^ 2$ contained in $X$  such that $Y_E\in |\mathcal{O}_E (1)|$ and $Y'$ is just the image of $Y$ in $X'$.\end{definition}

\begin{theorem}\label{6.2} Let $X,Y$ be as above and assume that $K_Y$ is not nef. Then one of the following holds.

\begin{itemize}
\item[{\rm (i)}] $\varrho (X) =1$, $X$ is Fano, of index $\ge 2$ and either:
\begin{itemize}
\item[{\rm (a)}] $X\cong \mathbb{P}^3$, $Y \in |\mathcal{O}_{\mathbb{P}^3}(a)|$, $a= 1,2,3$; or
\item[{\rm (b)}] $X\cong \sQ^3$, $Y\in |\mathcal{O}_{\sQ^3}(a)|$, $a=1,2$; or
\item[{\rm (c)}] $X$ is a del Pezzo three fold, $Y \in |\mathcal{O}_X(1)|$, cf. {\rm \cite{FuBook}} or {\rm \cite{IP}} for a complete 
list.\end{itemize}
\item[{\rm (ii)}] $X$ is a linear $\pn 2$-bundle over a curve and, for each fiber $F$, either $Y_F \in |\sO_{\mathbb{P}^2}(1)|$ or $Y_F\in|\sO_{\mathbb{P}^2}(2)|$;
\item[{\rm (iii)}] $X$ admits a contraction of an extremal ray, $\varphi :X\to W$, such that $W$ is a (smooth) curve, we have $F\cong\sQ ^2$ for a 
general fiber of $\varphi$ and $Y_F\in |\sO_{\sQ ^2} (1)|$ (we call $\varphi$ a {\em quadric fibration});
\item[{\rm (iv)}] $X$ is a linear $\pn 1$-bundle over a surface and $Y$ is a rational section;
\item[{\rm (v)}] A reduction $(X',Y')$ of $(X,Y)$ exists.
\end{itemize}
\end{theorem}
\begin{proof}
Since $K_X+Y$ is not nef, there exists an extremal ray $R= \reals_+ [C]$ of $X$ such that $(K_X+Y) \cdot C
<0$. Consider the length $\ell(R)$ of $R$ and observe that we have $\ell(R)\geq 2$. Let $\varphi ={\rm cont}_R :X\to W$ be the contraction of $R$ and let $F$ be a general fiber of $\varphi$. If $\dim W =0$, 
we fall in case (i). So, from now on, we may assume $\dim W >0$.  If $\ell(R) =4$, by Theorem~\ref{ExRayIW},
$\dim W=0$. If $\ell (R)=3$, by Theorem~\ref{ExRayIW}, $W$ is a curve. Moreover, $Y\cdot C=1$ or $2$. By Corollary~\ref{Co}, $F\cong \pn 2$. If $Y\cdot C =1$, we get case (ii), $Y_F \in |\sO_F(1)|$ by Theorem~\ref{FB}. Assume now that $Y\cdot C=2$ (and $\ell(R)=3$). Let $L:= K_X+2Y$ and let $F_0$ be an arbitrary fiber of $\varphi$. Remark that $L\cdot R > 0$, therefore $L\cdot C_0> 0$ for any curve $C_0\subset F_0$. We have that $(L_{F_0})^2=1$ and $L_{F_0}$ is ample by the Nakai--Moishezon criterion. By Theorem~\ref{KOF}, $\varphi$ makes $X$ a $\pn 2$-bundle and $Y_F\in|\sO_{\mathbb{P}^2}(2)|$. Thus, when $\ell(R)=3$ and $W$ is a curve, we get case (ii). Next, suppose that $\ell(R)=2$, so $Y\cdot C=1$. If $W$ is a curve we get $K_F +2Y_F\sim  0$ and we deduce from Remark~\ref{Re} that $F\cong \sQ^2$, leading to case (iii). If $W$ is a surface, we get case (iv). Indeed, $\varphi$ is generically a $\pn 1$-bundle, $Y$ being a rational section. So it is enough to see that all fibers of $\varphi$ are one-dimensional. Let $S$ be an irreducible surface contracted by $\varphi$ and let $D:=S\cap Y$. We obtain $S^2\cdot Y=(D_Y)^2<0$ and $S^2\cdot Y=D\cdot S=0$ since $S$ is contracted by $\varphi$. This contradiction shows that $\varphi$ is a $\pn 1$-bundle. Finally, assume that $\varphi$ is birational.
For such a ray, it follows from Theorem~\ref{ExRayIW} that $E$, the locus of $R$, is an 
irreducible surface, contracted to a point. Moreover, $E\cdot C:= -c<0$ since $R$ is not nef. We get $K_E+(c+2)Y_E \sim  0$;  as
above, using suitable vanishings, (see \cite{IoGen} for details) we deduce that $E\cong \pn 2$, $E_E \in |\sO_{\pn 2} (-1)|$ and $Y_E \in |\sO_{\pn 2}(1)|$. This leads to the 
reduction from case (v).
\end{proof}

\begin{corollary}\label{6.3}{\rm (\cite{BadAmp1, BadAmp2, BadAmp3})} Let $(X,Y)$ be as above and assume that $p:Y\to B$ is a $\pn 1$-bundle. Then $p$ 
extends to a linear $\pn 2$-bundle $\overline p: X\to B$, unless either $X\cong\pn 3$, $Y\in |\sO_{\pn 3}(2)|$, or $X\cong \sQ^3$, $Y\in 
|\sO_{\sQ^3}(1)|$, or $Y \cong \pn 1 \times \pn 1$, $p$ is one of the projections and the other projection extends.
\end{corollary}
\begin{proof}
Assume first that $\varrho (X) =1$. The conclusion follows by looking at the list in Theorem~\ref{6.2}(i), using the classification of del Pezzo 
three folds. Next, suppose $\varrho(X) >1$. By Lefschetz's theorem, 
we 
get an isomorphism ${\rm Num}(X)\cong {\rm Num}(Y)$. Then Corollary~\ref{newcor} applies to give that some extremal ray of $Y$ 
extends to $X$. If the genus $g(B) >0$, such a ray is unique and its contraction, $p$, extends. Apply Theorem~\ref{FB} to conclude. Assume 
$g(B)=0$. Unless either $Y\cong \pn 1\times \pn 1$  or $Y\cong \mathbb{F}_1$, $Y$ has only one extremal ray, so the previous argument applies. To conclude, 
we only have to examine the case when the contraction of the $(-1)$ curve of $\mathbb{F}_1$, say $\pi :\mathbb{F}_1\to \pn 2$, extends to a morphism $\overline \pi: X\to W$.

{\it Case} 1. $\overline\pi: X\to W$ is the contraction of a $(-1)$ plane $E$. The diagram 
\[\begin{xy}\xymatrix{
Y \ar[d]_{\pi}\ar@{^{(}->}[r] & X\ar[d]^{\overline\pi} \\
\pn 2\cong Z\ar@{^{(}->}[r]& W}\end{xy}\]
shows that $W \cong \pn 3$, $Z \in |\sO_{\pn 3}(1)|$  (see Theorem~\ref{Pn}). Let $L:= \overline{\pi}{\,}^*\sO_{\pn 3}(1)$. 
We get $Y^3=(L-E)^3=0$, contradicting ampleness of $Y$.

{\it Case} 2. $\overline\pi: X\to \pn 2$ is a $\pn 1$-bundle, $Y$ is a rational section (see the argument from the proof of Theorem~\ref{6.2}, case (iv)). 
Let $f$ be a fiber of $p$ and let $C_0$ be the ($-1$) curve contracted
by $\pi$. Write $Y_Y\approx aC_0 + bf$, for some $a>0$. From $Y\cdot C_0=1$ it follows $b=a+1$. Let $L:= \overline{\pi}{}^*(l)$, 
$l\subset\pn 2$ a line. Since $\Pic X \cong \Pic Y$, there is some $F\in \Pic X$ such that $\sO_Y (F) \cong \sO_Y(f)$.
We find easily that $Y\approx F+aL$. Now consider the exact sequence
\[0\to -aL \to \sO_X(F)\to \sO_Y(F)\to 0.\]
We have $H^1(X,-aL)=H^1(\pn 2 , -al)=0$. Therefore, using also the ampleness of $Y$, it follows that the linear system $|F|$ gives a morphism $\overline p: X\to \pn 1$ which extends $p$. Clearly, $\overline p $ is a $\pn 2$-bundle and, in fact, $X\cong \pn 1\times \pn2$. 
 \end{proof}

A classification of all cases when $Y$ is birationally ruled also follows from Theorem~\ref{6.2}.

\begin{corollary} {\rm (cf. \cite{Somin, SoVar})} Let $(X,Y)$ be as above and assume that $Y$ is not birationally ruled. Then, either $X$ is a 
(linear) $\pn 1$-bundle and $Y$ is a rational section, or there is a reduction $(X_0,Y_0)$ such that $K_{Y_0}$ is nef. 
\end{corollary}
\begin{proof}
Looking over the cases {\rm (i)}--{\rm (v)} in Theorem~\ref{6.2} and using the hypothesis that $Y$ is not ruled, we see that only cases {\rm (iv)} 
and {\rm (v)} are possible.
\end{proof}

\begin{corollary}
{\rm (\cite{IoGen})} Let $(X,Y)$ be as above. Assume that $X$ is not a $\pn 1$-bundle, $Y$ being a rational section. 
\begin{itemize}
\item[{\rm (i)}] If $\grk(Y)=0$, there is a reduction $(X_0,Y_0)$ of $(X,Y)$ such that $Y_0$ is a $K3$ surface
and $X_0$ is Fano.

\item[{\rm (ii)}] If $\grk(Y)=1$, there is a reduction $(X_0, Y_0)$ of $(X,Y)$ such that $X_0$ fibers over a curve, with general fiber a del Pezzo
surface.
\end{itemize}
\end{corollary}

\begin{proof}
Use the preceding corollary and Proposition~\ref{IoCY}. \end{proof} 

\section{Extending $\pn 1$-bundles}\label{conicbd}\addtocounter{subsection}{1}\setcounter{theorem}{0}

We start with the following proposition.

\begin{prop}\label{conicbdlethm} Let  $X$ be a projective manifold of dimension $n\geq 4$. Let $Y$ be a smooth ample divisor on $X$. Assume that $Y$ is a conic fibration, with general fiber $f$.
Let $V$ be the family of rational curves induced by $f$ on $X$. Then the following conditions are equivalent:
\begin{enumerate}
\item[{\rm (i)}] $Y\cdot f=1$;
\item[{\rm (ii)}] $V$ is unsplit;
\item[{\rm (iii)}] $V$ is locally unsplit.
\end{enumerate} If  {\rm (i)}--{\rm (iii)} hold, then $p$ is a $\pn 1$-bundle which extends to $\overline{p}:X\to Z$. Moreover, $\dim Z=2$ and $\overline{p}$ is a $\pn 2$-bundle.
 Conversely, if $p$ is smooth and extends, conditions {\rm (i)}--{\rm (iii)} hold.
\end{prop}
\proof   Since ${\rm (i)} \Longrightarrow  {\rm (ii)}  \Longrightarrow {\rm  (iii)}$ are clear, it is enough to show ${\rm (iii)} \Longrightarrow {\rm (i)}$. Let $a:=Y\cdot f$, let $y\in Y$ be a fixed general point and consider the standard exact sequence
\[0\to{\textstyle\bigoplus\limits^{n-2}} \sO_f(-1)\to N_{f/X}(-1)\to\sO_f(a-1)\to 0.\] Since $h^1(N_{f/X}(-1))=0$,  general facts from deformation theory of rational curves say that $\dim V_y=h^0(N_{f/X}(-1))=a$, and hence $\dim \scf_y=a+1$, see \ref{families}. By semicontinuity, the same holds at a general point $x\in X$. Fix such a general point $x\in X$ and take another point $t\in {\rm Locus}(V_x)$. Since $V$ is locally unsplit, we know that each curve from $V_x$ is irreducible. By the non-breaking lemma, we thus conclude that there is a finite number of curves in $V_x$ passing through $t$. That is the projection  $q:\scf_x\to {\rm Locus}(V_x)$ is a finite map. Therefore $\dim {\rm Locus}(V_x)=a+1$. 
Thus we obtain $\dim Y\cap{ \rm Locus}(V_x)\geq a$. Assume by contradiction that $a\geq 2$. Then there exists a curve $C\subset Y\cap{ \rm Locus}(V_x)$ such that $p(C)$ is a curve in $Z$. In this case, a variant of the non-breaking lemma (see \cite[(1.14)]{Wlength} and also \cite[(1.4.5)]{BSW}) implies that the curve $C$ is numerically equivalent in $X$ to $\lambda f$ for some positive rational number $\lambda$.
Now, take a hyperplane section $H_Z$ of $Z$ and let $\sL\in \Pic X$ be the extension of $p^*(H_Z)$ on $X$ via the isomorphism $\Pic X\cong \Pic Y$. In particular, $\sL\cdot f=0$,  this leading to the numerical contradiction $ 0< \sL\cdot C=\lambda(\sL\cdot f)=0$.

If (i)--(iii) hold, $p$ extends to a $\pn 2$-bundle by Lemma~\ref{key1}. Moreover, $\dim Z=2$ by Lemma~\ref{easy}. Conversely,
if $p$ extends, we have (i) by Theorem~\ref{FB}.
\qed

We consider now the extension problem for $\pn 1$-bundles. For perspective we also recall the (much easier) case of $\pn d$-bundles, for $d\geq 2$.

The following major
conjecture on the topic \cite[Section 5.5]{Book} describes all known examples. 
We refer to \cite[Section 5.5]{Book} for the more general case when  $X$ is a local complete 
intersection. 
\begin{conjecture}\label{PdBdleConj} Let $L$ be an  ample line bundle on a
 projective manifold, $X$, of dimension $n\ge 3$. Assume that there is a
smooth
 $Y\in |L|$ such that $Y$ is a $\pn d$-bundle, $p:Y\to Z$, over a manifold, $Z$, of dimension $b$.  Then 
$d\ge b-1$ and   $(X,L)\cong ({\mathbb P}(\sE),\xi_{\mathbb P})$ for an ample vector bundle, $\sE$, 
on $Z$ with $p$ equal to the
restriction to $Y$ of the induced projection ${\mathbb P}(\sE)\to Z$, except if either:
\begin{enumerate}
\item[{\rm (i)}]  $(X,L)\cong (\pn 3, \sO_{\pn 3}(2))$; or
\item[{\rm (ii)}]  $(X,L)\cong (\sQ^3, \sO_{\sQ^ 3}(1))$; or
\item[{\rm (iii)}]  $Y\cong \pn 1\times \pn{n-2}$, $p$ is the product projection onto the second factor, $(X,L)\cong 
(\proj\sE,\xi_{\pn{}})$ for an ample vector bundle, $\sE$, 
on $\pn 1$ with the product projection of $Y$ onto the first factor equal to  the induced projection
$\proj\sE\to \pn 1$.
\end{enumerate}
\end{conjecture}

\noindent{{\em Note.}} The inequality $d\geq b-1$ is a necessary condition for $p:Y\to Z$ to extend by
Lemma~\ref{easy}. 

\smallskip

The conjecture has been shown except when $d=1$, $b\ge 3$, and the base $Z$ does not map finite-to-one into its Albanese variety.  
The case when either $d\ge 2$ or $Z$ is a submanifold of an abelian variety 
follows from  Sommese's
extension theorems   \cite{So1} (see also Fujita \cite{FuLef}). This argument works also in the case when $Z$ maps 
finite-to-one into its Albanese variety (see  \cite[(5.2.3)]{Book}).

\begin{theorem}\label{d2} {\rm (Sommese)} Conjecture~\ref{PdBdleConj} is true for $d\geq
2$.\end{theorem}
\proof  Since the result is trivial if $Z$ is a point we can assume without loss of generality that
 $\dim Z\ge 1$ and thus that $n\ge d+2\ge 4$.
 From  Theorem~\ref{ExtensionThm} we know  that 
$p : Y\to Z$ extends to a morphism, $\overline{p} : 
X\to Z$.  The result follows from Theorem \ref{FB}(ii).
\qed 

The conjecture is also known when $d=1$ and $b\le 2$. If $b=1$, the result is due to B\u adescu \cite{BadAmp1, BadAmp2, BadAmp3};
we have seen a proof in Corollary~\ref{6.3}. The case $b=2$ is due to the contribution of several authors: Fania and Sommese \cite{FS},
Fania, Sato and Sommese \cite{FSS}, Sato and Spindler \cite{SatSpin1, SatSpin2} and also \cite {Sat1, SaYi}. Below we propose a shorter proof.
The basic ideas are those in \cite{FSS} and \cite{SaYi}, but we do not use the difficult papers \cite{FS} and \cite{Sat1}.

\begin{theorem}\label{mainPd}
Conjecture~\ref{PdBdleConj} is true when $d=1$ and $b=2$.
\end{theorem}
\begin{proof}
Assume that $p$ does not extend. 

{\it  Step} 1. $Z$ is ruled. From Corollary~\ref{newcor} it follows that $Y$ has some extremal ray, say $R$, which extends to an extremal ray 
$\overline R$ on $X$. Using \cite{Mo1}, we consider the possible type of $R$. If $R$ is nef, $Z$ is covered by rational curves
in the fibers of ${\rm cont}_R$, so it is ruled. If $R$ is not nef, $E$, the locus of $R$, covers $Z$ (so again $Z$ is ruled), unless 
$E \cong \pn 1\times \pn 1$ and $p(E):=C$ is a curve. Standard computations (cf. also Poposition~\ref{Palbup2}) show that 
$C$ is a $(-1)$ curve and we may construct a commutative diagram 
\[\begin{xy}\xymatrix{Z\ar[d]_{{\rm cont}_C} &\ar[l]_p Y \ar[d]^{{\rm cont}_R}  \ar@{^{(}->}[r]& X\ar[d]^{{\rm cont}_{\overline R}}\\
Z'&\ar[l]^{p'} Y'  \ar @{^{(}->}[r] &X'}\end{xy}\]
where $Y'$ is ample on $X'$ (cf. \cite{FuLef}) and $p'$  is again a $\pn 1$-bundle. So, after finitely many steps, we conclude that 
$Z$ is ruled.

{\it Step} 2. $Z\cong \pn 2$. Assume the contrary. As $Z$ is ruled, there is a morphism $\varphi:Z\to B$ which is generically a $\pn 1$-bundle.
Apply Theorem~\ref{ExtensionThm} to the map $\pi: =\varphi\circ p$ to get an extension $\overline{\pi}:X\to B$. Next we use 
Corollary~\ref{6.3} fiberwise. Let $F\cong \mathbb{F}_e$, $\overline F$ be the general fibers of $\pi, \overline{\pi}$ respectively.
Denote by $f$, $C_0$ a fiber and a minimal section of $F$, respectively.
Note that the classes of $f$ and $C_0$ are not proportional in $N_1(Y)$. Then the diagram
\[\begin{xy}\xymatrix{N_1(F)\ar[r]\ar[d]&N_1(\overline F)\ar[d]\\N_1(Y)\ar[r]^\sim &N_1(X)}\end{xy}
 \]  
 shows that $\dim_{\reals}N_1(\overline F)\geq 2$. So, from Corollary~\ref{6.3} we infer that either $Y\cdot f=1$, or $F\cong \pn 1\times \pn 1$. In the first case, $p$ extends by Lemma~\ref{key1}. So we may assume that $F\cong \pn 1\times \pn 1$ from now on.
 
 We refer to \cite[pp.\ 7--11]{FSS} for details concerning the next few arguments. Using standard properties of Hilbert schemes, one shows:
 \begin{itemize}
 \item[(a)] $\varphi $ is a $\pn 1$-bundle;
 
 \item[(b)] any fiber of $\pi$ is isomorphic to $\pn 1\times \pn 1$;
 
 \item[(c)] the family of curves on $Y$ determined by minimal sections of the map $p_{|F}:F\to p(F)$ yields another $\pn 1$-bundle $p':Y \to Z'$. We deduce a cartesian diagram 
 \[\begin{xy}\xymatrix{&\ar[dl]_p Y\ar[dr]^{p'}&\\Z\ar[dr]_\varphi &&Z'\ar[dl]^{\varphi'}\\&B&}\end{xy}
 \]
 where $\varphi' $ is also a $\pn 1$-bundle.
 \end{itemize}
 
 Next, from the  above construction, Corollary~\ref{6.3} and Theorem~\ref{FB}, we find that $p'$ extends to a linear $\pn 2$-bundle $\overline p{}':X\to Z'$. 
 \begin{itemize}
 \item[(d)] This yields an exact sequence \[
 0\to \sO_{Z'} \to \sF\to \sG\to 0,\]
 where $\sF, \sG$ are ample vector bundles on $Z'$. 
 If $Z'=\pn{}(\sE')$, one finds $\sG\cong \varphi'{}^*(\sE) \otimes \xi^a$, where $\sE$ is a rank $2$ vector bundle on $B$, $\xi=\xi_{Z'}$ and 
 $a>0$. Moreover, from (c) it follows that $Z\cong \pn{}(\sE)$.
 \item[(e)] If we assume $\sE$ unstable, it is now standard to see that the exact sequence  from (d) splits.
 This is a contradiction, since $\sF$ is ample.
   \end{itemize}
   
   Finally, see \cite{SatSpin1} for a proof that the case $\sE$ stable also leads to a contradiction.
   
   {\it Step} 3. Conclusion. We know that $Z\cong \pn 2$, so Proposition~\ref{FSPb} below applies to give the result.
   \end{proof}
   
Let us also explicitly note that in the  relevant case $d=1$, $b\geq 3$ (by Lemma~\ref{easy} the bundle $p:Y\to Z$
does not extend in this case) Conjecture~\ref{PdBdleConj} is equivalent to saying that
\begin{itemize} \item {\em A $\pn 1$-bundle $Y$, $p:Y\to Z$, over a manifold $Z$ cannot be an ample divisor in an
$n$-dimensional manifold $X$ unless $Y\cong\pn 1\times\pn {n-2}$, $Z\cong \pn {n-2}$, and $X$  is a $\pn
{n-1}$-bundle over
$\pn 1$ whose restriction to $Y$ is the projection $Y\to \pn 1$ on the first factor}.\end{itemize}


The following proposition from \cite{FSS} ensures that,  to prove Conjecture~\ref{PdBdleConj}, it is enough to show that $Z\cong \pn {n-2}$, 
assuming that $p$ does not extend. 

\begin{prop}\label{FSPb} {\rm  (\cite[Section 2]{FSS})} Let $Y$ be a smooth ample divisor on an $n$-dimensional projective manifold $X$. Assume that 
 $p:Y\to Z$ is a  $\pn 1$-bundle over $Z=\pn b$, $b\geq 2$.  
If $p$ does not extend to $X$, then $Y\cong \pn 1 \times \pn b$ and $X$ is a $\pn {b+1}$-bundle over $\pn 1$.\end{prop} 
\begin{proof}
Let $F \in |p^* \sO_{\pn b}(1)|$ and let $L: = \sO_X(Y)$. Looking over the proof of Theorem~\ref{ExtensionThm} we see that, 
once we have the vanishings in \eqref{vanishing}, the argument works also in the case $\dim Y -\dim Z=1$. Therefore, the assumption that $p$ does not extend translates into: there exists some $t>0$ such that $H^1(Y, F-tL_Y)\not =0$. Using Serre duality, Kodaira vanishing, and the exact sequence
\[0\to K_Y+tL_Y-F \to K_Y+tL_Y \to K_F +tL_F-F\to0\]
it follows that
\begin{equation}\label{7.1}
H^{n-3}(F, K_F+tL_F-F)\not = 0\quad \mbox{for some }t>0.
\end{equation}
Iterating this construction, we may assume that $b=2$; in this case $F=\mathbb{F}_e$, for some $e\geq 0$. We write $Y=\pn{}(\sE)$ for 
some rank $2$ vector bundle on $\pn b$. We also may assume that, if $l \subset \pn b$ is a line, we have $\sE_l\cong \sO_l\oplus \sO_l(-e)$.
We shall prove that $e=0$, so that $\sE$ is trivial (see \cite[Section~3.2]{OSS}) and the conclusion follows.

So, assume that $b=2$, $F=\mathbb{F}_e$ and write $L_F\sim  aC_0+bf$, using the notation from \cite[Chapter~V.2]{Hrt}. 
Since $L$ is ample, $a>0$ and $b>ae$. Hence, for $t>0$, $bt > aet$. Therefore, either $bt-1>aet$ and $tL_F -F\sim  atC_0+(bt-1)f$ 
is ample, or $bt=aet +1$ and $tL_F-F\sim  at(C_0+ef)$. Now, if $e>0$, $C_0 +ef$ is nef and big. Using Kawamata--Viehweg vanishing this contradicts \eqref{7.1}. So $e=0$ and we are done. 
\end{proof}

Further evidence for Conjecture~\ref{PdBdleConj} is given by the following result (see \cite{BFS2} for a proof).

\begin{prop}\label{Flemma} {\rm  (\cite{BFS2})} Let $p:Y\to Z$ be a  $\pn 1$-bundle over a smooth
projective three fold $Z$. Then $Y$ cannot be a very ample divisor in any projective
manifold $X$, unless $Z\cong\pn 3$ and $Y\cong\pn 1 \times\pn 3$.\end{prop}


\begin{rem*}\label{beautyfulmind} The following discussion gives some further support  to  Conjecture \ref{PdBdleConj},
in relation to the content of Section~4.  Let $p:Y\to Z$ be a smooth $\pn 1$-bundle with $Y$ ample divisor in a projective manifold $X$ of dimension $n\geq 4$. If $p$ does not extend to $X$, and we assume the hypothesis of Proposition~\ref{sexy} is fulfilled, then the following three conditions hold true.\begin{enumerate}
\item[{\rm (1)}] $\varrho(X)=\varrho(Y)=2$;
\item[{\rm (2)}] $Y$ is a Fano manifold;
\item[{\rm (3)}] $Z$ is a Fano manifold (and $\varrho(Z)=1$).
\end{enumerate}

Condition (1) directly follows from Proposition~\ref{sexy}.

To show (2), let $R_1$ be the extremal ray corresponding to the bundle projection $p$. By Corollary~\ref{newcor} we conclude that there exists an extremal ray $R_2$ on $Y$ which extends to $X$.  By (1), $ \overline{NE}(Y)=\langle R_1,R_2\rangle$ and $Y$ is a Fano manifold.

Since $Y$ is a Fano manifold and $p$ is smooth, $Z$ is also a Fano manifold, see \cite[p.\ 244]{KoBook}. In our special case, we can give the following alternative argument. Assume that  $Z$ is not a Fano manifold. Then $K_Z$ would be nef  (since  $\varrho(Z)=1$). Therefore Proposition~\ref{Lucian1} applies to say that $p$ extends; a contradiction. Whence  (3) holds.

We already observed that, in view of Proposition~\ref{FSPb}, to prove Conjecture~\ref{PdBdleConj} it would be enough to show that $Z\cong \pn {n-2}$. This is in agreement with the fact that assuming the hypothesis of Proposition~\ref{sexy} to hold, we get condition (3) above. 
For instance, when  $n=4$, the only Fano surface $Z$ with $\varrho(Z)=1$ is $\pn 2$, yielding a very short proof of Theorem~\ref{mainPd}. 
\end{rem*}

\section{Fano manifolds as ample divisors}\label{Fano}\addtocounter{subsection}{1}\setcounter{theorem}{0}

Throughout this section, let $X$ be a projective manifold of dimension $n\geq 4$, let $H$ be an ample
line bundle on $X$, and let $Y$ be a smooth divisor  in $|H|$. 
In this general setting, it is natural to ask  the following questions.

\begin{question}\label{Q1} If $Y$ is a Fano manifold, when is $X$  a Fano manifold?
\end{question}

\begin{question}\label{Q2} If $X$ is a Fano manifold, when do we have $\overline{NE}(X)\cong\overline{NE}(Y)$?
\end{question}

Question~\ref{Q2} has been solved in \cite{Ko} and \cite{Wcontr}, by using Theorem \ref{AW}, in the special case described in the following theorem.

\begin{theorem}\label{NE} {\rm (Koll\'ar, Wi\'sniewski)} Let $X$ be a Fano manifold of dimension $n\geq 4$ and index $i\geq 1$, $-K_X\cong iL$, for some ample line bundle $L$ on $X$. Assume that we are given a smooth member $Y\in |mL|$, for some integer $m$, $1\leq m\leq i$.  Then there is an isomorphism $\overline{NE}(Y)\cong \overline{NE}(X)$.\end{theorem} 
\proof  The result follows from Proposition~\ref{reph}, since $-(K_X+Y)\approx (i-m)L$ is nef.
\qed

Let us come back now to Question \ref{Q1}. This question is motivated by the problem of classifying polarized pairs $(X,Y)$ as above, when $Y$ is a Fano manifold of large index. Set
$-K_Y\approx iL_Y$, where $L_Y$ is an ample line bundle on $Y$.

The first cases to consider, cf. Corollary~\ref{Co} and Remark~\ref{Re}, are $i=\dim Y+1=n$, so $Y$ is a projective space, and $i=n-1$, so $Y$ is a quadric; for a solution, see \cite{So1} and also
\cite{BadHD}, where $X$ is only assumed to be normal. The case when $(Y,L_Y)$ is a classical del Pezzo variety, i.e., $i=n-2$ with
$L_Y$ very ample, has been completely worked out by Lanteri, Palleschi and Sommese \cite{LPS1,LPS2,LPS3}. In \cite{BF,BFS1} the next case when $(Y,L_Y)$ is a Mukai variety, i.e., $i=n-3$, is
considered. The results of \cite{BFS1} have been refined and strengthened in \cite{anoc} under the assumption
that $L_Y$ is merely ample, as a consequence of a comparing cones result which holds true in the range $i\geq
\frac{\dim Y}2$. In \cite{NO} the classification is extended to  the next case. 

We will work under the {\it extra assumption} that  the line bundle $L_Y$ is spanned. Note that this is
in fact the case when  $(Y,L_Y)$ is either a del Pezzo variety of degree at least two, or a Mukai variety. This follows from
Fujita's classification \cite{FuBook} of del Pezzo manifolds and from Mukai's classification, see \cite{Mu1} and \cite{Mella}.

We have the following result (compare with  \cite[(4.2)]{anoc}).

\begin{theorem}\label{BI} Let $X$ be a projective manifold of dimension $n\geq 4$, let $H$ be an ample
line bundle on $X$, and let $Y$ be an effective divisor  in $|H|$. Assume that $Y$ is a Fano manifold of index
$i\geq 3$, $-K_Y\approx iL_Y$. Further assume  that $L_Y$ is spanned. Then either
\begin{enumerate}
\item[{\rm  (i)}] There exists an extremal ray $R$ on $X$ of length $\ell(R)\geq 2i+1$; or
\item[{\rm (ii)}] $Y$ contains an extremal
ray of length $\geq 3i$; or
\item[{\rm (iii)}] $X$ is a Fano manifold and $\overline{NE}(X)\cong \overline{NE}(Y)$.\end{enumerate}\end{theorem}

\proof By the Lefschetz theorem, there exists a unique line bundle $L$ on $X$ such that $L_{|Y}=L_Y$. 

First, suppose that $K_Y+iH_Y$ is not nef. Then, by the cone theorem, there exists an extremal ray
$R_Y=\reals_+[C]$ on $Y$ such that 
\[(K_Y+iH_Y)\cdot C<0.\]
We can also assume that $C$ is a minimal curve,
that is $\ell(R_Y)=- K_Y\cdot C$. Therefore  
$\ell(R_Y)>i(H\cdot C)$.

Set $a:=L\cdot C$, $d:=H\cdot C$. Since $-K_Y\approx iL_Y$, the last inequality yields $a>d\geq 1$. Then
$a\geq 2$ with equality only if $d=1$. Thus $\ell(R_Y)=ia\geq 2i$, so that either $\ell(R_Y)\geq 2i+1$,
or $a=2$ and hence $d=H\cdot C=1$.

From $-K_Y\approx iL_Y$ we get $\ell(R_Y)=i(L_Y\cdot C)$, that is, $i$ divides $\ell(R_Y)$. Hence in the first
case above it must be $\ell(R_Y)\geq 3i$, as in the case (ii) of the statement.

Thus we can assume that $a=2$, $H\cdot C=1$ and the adjunction formula gives
$(K_X+H)\cdot C=K_Y\cdot C=-\ell(R_Y)$, or $-K_X\cdot C=\ell(R_Y)+1=2i+1$.
Therefore $(K_X+2iH)\cdot C<0$, i.e., $K_X+2iH$ is not nef. Then there exists a rational curve $\gamma$
generating a ray $R=\reals_+[\gamma]$ on $X$ such that $(K_X+2iH)\cdot\gamma<0$. Since  we can assume
$\ell(R)=-K_X\cdot \gamma$, it follows that $\ell(R)>2i$ and we are in case (i) of the statement.

Assume now that $K_Y+iH_Y$ is nef. Then by the ascent of nefness (see the proof of Theorem~\ref{six}) we infer that
$K_X+(i+1)H$ is nef and hence by Kawamata--Reid--Shokurov base point free theorem we conclude that $m(K_X+(i+1)H)$ is
spanned for $m\gg 0$.

We proceed by cases, according to the Iitaka dimension of $K_X+(i+1)H$. Let $\psi:X\to W$ be the map with
normal image and connected fibers associated to $|m(K_X+(i+1)H)|$ for $m\gg 0$.

If $\grk(K_X+(i+1)H)=0$, then $K_X+(i+1)H\approx 0$, so that $X$ is a Fano manifold. From $iL_Y\approx
-K_Y\approx iH_Y$ we conclude by Lefschetz that $H\approx L$. Therefore Theorem \ref{NE} applies to
give $\overline{NE}(X)=\overline{NE}(Y)$.

Assume $\grk(K_X+(i+1)H)=1$. Since $Y$ is ample, the restriction $\psi_Y$ of $\psi$ to $Y$ maps onto $W$, so
that $W\cong \pn 1$ since $Y$ is a Fano manifold. Note that $\psi_Y$ is not the constant map by ampleness of
$Y$. Recalling that $\overline{NE}(Y)$ is polyhedral, we conclude that there exists an extremal ray $R$ on $Y$
which is not contracted by $\psi_Y$. Let $\vphi:Y\to Z$ be the contraction of $R$. We claim that all fibers of
$\vphi$ are of dimension $\leq 1$. Otherwise, let $\Delta$ be a fiber  of dimension $\geq 2$. Any fiber
$F$ of
$\psi_Y$ is a divisor on $Y$. Then we can find a curve $C\subset \Delta\cap F$. Therefore $C$ generates $R$
and $\dim\psi_Y(C)=0$, contradicting the fact that $R$ is not contracted by $\psi_Y$.

Thus by Theorem \ref{AW} we know that either $\vphi$ is a blowing-up of a smooth codimension two subvariety of
$Z$ and $-K_Y\cdot C=1$; or $\vphi$ is a conic fibration and $-K_Y\cdot C\le 2$. In each case, the equality
$-K_Y\cdot C=i(L_Y\cdot C)$ contradicts the assumption that $i\geq 3$.

 Assume now $\grk(K_X+(i+1)H)\geq 2$. We follow here the argument from \cite{BFS1}. From $(K_X+H)_Y\approx
 K_Y\approx -iL_Y$ we get by Lefschetz
 \begin{equation}\label{basic}
  K_X+H+iL\approx 0,\end{equation} that is $K_X+(i+1)H\approx i(H-L)$. 
Thus we conclude that  $\grk(H-L)\geq
 2$ and that $m(H-L)$ is spanned for $m\gg 0$. Therefore the Mumford vanishing theorem \cite{Mum} (see also \cite[(7.65)]{SS}) applies
to give
\begin{equation}\label{SSvan} H^1(X, L-H)=0.\end{equation}
Now consider the exact sequence
\[0\to L-H\to L\to L_Y\to 0.\] Since $L_Y$ is spanned on $Y$,  by
(\ref{SSvan}) we see that sections of $H^0(Y, L_Y)$ lift to span $L$ in a neighborhood of
$Y$; but since $Y$ is ample we conclude that $L$ is spanned off a finite set
 of points. Hence $L$ is nef and therefore $-K_X$ is ample by (\ref{basic}), i.e., $X$ is a Fano manifold.

We conclude that either 
$\overline{NE}(X)\cong \overline{NE}(Y)$ and hence we are done, or there exists an extremal ray $R=\mathbb{R}_+[C]$ on $X$,
$R\subset \overline{NE}(X)\setminus \overline{NE}(Y)$ such that every fiber of the contraction $\vphi:X\to Z$ of $R$ has dimension at most one. Then
Theorem~\ref{AW} applies again  to say that either:
\begin{enumerate}
\item[(1)] $\vphi$ is a blowing-up along a smooth codimension two center $B$ and $K_X\cdot C=-1$, or
\item[(2)] $\vphi$ is a conic fibration and  either  $K_X\cdot C=-2$ or $K_X\cdot C=-1$. \end{enumerate}

In case (1), from (\ref{basic}) and $K_X\cdot C=-1$ we get 
$1=Y\cdot C+i(L\cdot C)$.     Since $Y\cdot C>0$ and $L\cdot C\geq 0$ it must be $Y\cdot C=1$, $L\cdot C=0$.  Note that $(K_X+Y)\cdot C=0$ and apply Proposition \ref{KollarProp} to contradict our present assumption  that $R\not\subset\overline{NE}(Y)$.

Let us consider case (2). If  $K_X\cdot C=-2$ we have $2=Y\cdot C+i(L\cdot C)$,  giving $Y\cdot C=2$, $L\cdot C=0$.  If  $K_X\cdot C=-1$ we get $1=Y\cdot C+i(L\cdot C)$,  giving $Y\cdot C=1$, $L\cdot C=0$. In both cases, Proposition~\ref{KollarProp} applies again to give the same contradiction as above. \qed

\begin{corollary}\label{anocor} {\rm (\cite[(4.2)]{anoc})} Let $X$ be an $n$-dimensional projective manifold. Let $H$ be an ample
line bundle on $X$ and let $Y$ be a divisor  in $|H|$. Assume that $Y$ is a Fano manifold of index  
$i\geq\frac{\dim Y}2\geq 3$  {\rm (}hence $n\geq 7${\rm )}, $-K_Y\approx iL_Y$ and  $L_Y$ is spanned. Then $X$ is a 
Fano manifold and $\overline{NE}(X)\cong\overline{NE}(Y)$. 
\end{corollary}

\proof  By the proof of  Theorem~\ref{BI}, either we are done, or $Y$ contains an extremal ray $R_Y=\reals_+[C]$ such that either $\ell(R_Y)\geq 3i$, or   $\ell(R_Y)\geq 2i$ and $H\cdot C=1$. In the first case we get  the numerical contradiction
\[\ell(R_Y)\geq \frac{3}{2}(n-1)\geq n+1=\dim Y+2.\]

Thus we may assume $\ell(R_Y)\geq 2i$ and $H\cdot C=1$.  Therefore,  if $\Delta$ is a positive dimensional fiber of the contraction $p={\rm cont}_{R_Y}:Y\to W$, we have $\dim\Delta\geq\ell(R_Y)-1\geq \dim Y-1$ (see Theorem~\ref{ExRayIW}).  Thus, either $\Delta=Y$ or  $\ell(R_Y)=\dim Y=n-1$.

In the first case, the contraction $p$ is the constant map, so that $\Pic Y\cong\Pic X\cong \zed$ and the conclusion is clear.

In the latter case, since $\ell (R_Y)=\dim Y$, we know from Theorem~\ref{ExRayIW} that $\dim W\leq 1$. If $\dim W=0$ we conclude as above. Assume that $W$
is a curve and let $F$ be a general fiber of $p$. Since $Y$ is a Fano manifold, $W\cong \pn 1$. Moreover, since $H\cdot C=1$, we get $K_F+(n-1)H_F\approx 0$. 
Corollary~\ref{Co} applies to give $F\cong \pn {n-2}$, $H\in |\sO_F(1)|$. Therefore, by Theorem~\ref{FB}(i), $Y\cong \pn{}(\sE)$ for some vector bundle $\sE$ on $\pn 1$. Using Lemma~\ref{wk} and the assumption $i\geq 3$, we see that this case does not occur.   
\qed

\section{Ascent properties}\label{Ascent}\addtocounter{subsection}{1}\setcounter{theorem}{0}

Let $X$ be a projective $n$-dimensional  manifold and let $Y\subset X$ be a smooth ample divisor.
Here is a list of general facts concerning ascent properties  from  $Y$   to $X$. E.g., 

$\bullet$ $K_Y$ not ample  $\Longrightarrow$ $K_X$ not nef. It immediately follows from the adjunction formula.

$\bullet$ $\grk(Y)<\dim Y$  $\Longrightarrow$ $\grk(X)=-\infty$. Here is the argument from \cite[Proposition~5]{IoGen}. Assume by contradiction that 
$|mK_X|\neq\emptyset$ for some $m>0$ and let $E\in |mK_X|$. Write $E=aY+E'$, with $a\geq 0$ and $Y
\not\subset{\rm Supp}(E')$. Since $\lambda Y$ is very ample for $\lambda\gg 0$, we can find $\lambda\gg 0$ and
$D\in |\lambda E|$ such that $Y \not\subset{\rm Supp}(D)$. Since we have $(D+\lambda mY)_{|Y}\approx \lambda
mK_Y$, it follows that $|\lambda mK_Y|$ is very ample outside of $Y\cap {\rm Supp}(D)$, so that $\grk(Y)=n-1=\dim
Y$. This contradiction proves the assertion.

Uniruled manifolds are birationally Fano fibrations. This fact follows from Campana's construction (see e.g., the Preface and Chapters $3$, $4$ of Debarre's text \cite{Debarre}). Many results in our paper are concerned with the case in which $Y$ carries a special Fano fibration structure. 

$\bullet$ $Y$ uniruled $\Longrightarrow$ $X$ uniruled. It immediately follows from the uniruledness criterion \cite[II, Section 3, IV, (1.9)]{KoBook}. Saying that $Y$  is
uniruled means that there is a morphism $f:\pn 1\to Y$ such that $f^*T_Y$ is spanned.
Consider the tangent bundle sequence
\[0\to T_Y\to T_{X|Y}\to \sO_Y(Y)\to 0.\]
  Let $f':\pn 1\to X$ be the induced morphism to $X$. By pulling back to
$\pn 1$, we get the exact sequence
\[0\to f^*T_Y\to f^*(T_{X|Y})={f'}^*T_X\to f^*\sO_Y(Y)\to 0.\]
Since both $f^*\sO_Y(Y)$ and $f^*T_Y$ are nef, we conclude that ${f'}^*T_X$ is  nef and hence spanned; this is 
equivalent to say  that $X$ is uniruled.

$\bullet$ $Y$ rationally connected  $\Longrightarrow$ $X$ rationally connected. Saying that $Y$ is rationally connected is equivalent to the existence of  a curve $C\cong \pn 1\subset Y$  with ample
normal bundle $N_{C/Y}$ (see e.g., \cite{KoBook}). Therefore the exact sequence of normal bundles
\[0\to N_{C/Y}\to N_{C/X}\to \sO_C(Y)\to 0\] and the ampleness of $Y$ give the ampleness of $ N_{C/X}$.

$\bullet$ $Y$ unirational $\Longrightarrow$ $X$ unirational? This is a hard question and no answer is known. It is interesting to point out that, since unirationality implies rational connectedness, to find examples of $Y$ unirational with $X$ not unirational would give  examples of rationally connected manifolds $X$  which are not unirational. Quoting  Koll\'ar \cite[Section 7, Problem 55]{KoVar}, the latter is ``one of the most vexing open problems'' in the theory.

$\bullet$ In general, $Y$ rational does not imply that $X$ is rational. We present below a few results about this problem. In particular, we  obtain a proof of the following classical statement  (\cite[Chapter IV]{Roth}): {\em for a  very ample  
smooth divisor  on a three fold $X$, the ascent of rationality holds true with the only exception when $X$ is the  cubic hypersurface of $\pn 4$}. The case of the cubic hypersurface is indeed an exception, see \cite{ClGr}.

\begin{theorem}\label{rat} {\rm (cf. also \cite{CaFl})} Let $L$ be an  ample line bundle on a
 smooth projective three fold $X$. Assume that there is a
smooth $Y\in |L|$ such that $Y$ is rational. Then $X$ is rational unless either:
\begin{enumerate}
\item[{\rm (i)}] $L^3=1$ and $(X,L)$ is a weighted hypersurface of degree $6$ in the weighted projective space
${\mathbb P }(3,2,1,1,1)$, $-K_X\approx 2L$; or
\item[{\rm (ii)}] $L^3=2$ and $(X,L)$ is the double covering of $\pn 3$ branched along a smooth surface of degree $4$,
$-K_X\approx 2L$ and $L$ is the pull-back of $\sO_{\pn 3}(1)$; or
\item[{\rm (iii)}] $X$ is the hypercubic in $\pn 4$ and $L\approx \sO_X(1)$.\end{enumerate}\end{theorem}

\proof Since $Y$ is a rational surface, $K_Y$ is not nef.  We follow the cases (i)--(v) from Theorem~\ref{6.2}. In case (i), we apply the well-known 
classification of Fano three folds of index  $\geq 2$ (and $\varrho(X)=1$), see \cite{IP}.
We either get one of the exceptional cases in the statement, or $X$ is the complete intersection of two quadrics 
in $\pn 5$, or $X \subset \pn 6$ is a linear section of the Grassmannian of lines in $\pn 4$, embedded in $\pn 9$ by the Pl\"ucker embedding. 
In the last two cases $X$ is rational (see e.g.,\ \cite{IP}). A simple argument is given in Example~\ref{sad} below.

Assume now that we are in case (ii) or (iii) from Theorem~\ref{6.2}. For such a fibration the base
curve is $\pn 1$ and the general fiber is rational. Moreover, a section exists by Tsen's theorem, see e.g.,\ \cite[IV.6]{KoBook}.
So $X$ is rational, too. 

If we are in case (iv), the base surface is birational to $Y$, so it is rational. We conclude that $X$ is rational.
Finally, case (v) leads to one of the previously discussed situations. 
\qed

The following  result, contained in  \cite[Theorem~1.3]{IoNa}, concerns   the ascent 
of rationality  from a suitable rational submanifold. The proof relies on Hironaka's desingularization theory \cite{Hir} and on basic properties of rationally connected manifolds \cite{KoMiMo}.

\begin{theorem}\label{trei-sase} {\rm (\cite[Theorem~1.3]{IoNa})}
Let $X$ be a projective variety and $|D|$ a complete linear system of Cartier divisors on it. Let $D_1,\ldots,D_s \in |D|$ and put $W_i := D_1\cap \cdots \cap D_i$ for $1\leq i\leq s$. Assume that $W_i$ is smooth, irreducible of dimension $n-i$, for all $i$. Assume moreover that there is a divisor $E$ on $W := W_s$ and a linear system $\Lambda \subset |E|$ such that:\begin{enumerate}
\item[{\rm (i)}] $\varphi_{\Lambda} : W \dasharrow \pn {n-s}$ is birational, and
\item[{\rm (ii)}] $|D_W - E|\neq \emptyset $.\end{enumerate}
 Then $X$ is rational.
\end{theorem}
\proof  We proceed by induction on $s$. Let us  explain the case $s=1$, the general case being
  completely similar. So, let $W \in |D|$ be a smooth, irreducible Cartier
  divisor such that $\varphi_{\Lambda} : W \dasharrow \pn {n-1}$ is birational for
  $\Lambda \subset |E|$, $E \in {\rm Div}(W)$ and $|D_W - E|\neq
  \emptyset$. Replacing $X$ by its desingularization, we may assume that $X$
  is smooth. As $W$ is rational, it is rationally connected, so we may find
  some smooth rational curve $C \subset W$ with $N_{C/W}$ ample. We have $C\cdot E > 0$ and from (ii) we deduce $C \cdot D > 0$. From the exact sequence of normal bundles we get that $N_{C/X}$ is ample, so $X$ is rationally connected. In particular, $H^1(X, \sO _{X}) = 0$. 
  
  The exact sequence 
  \[0 \to \sO_X \to \sO_X (D) \to \sO_W (D) \to 0,\]
  shows that $\dim |D| = \dim |D_W| + 1 \geq \dim |E| + 1 \geq n$.
  
  We may choose a pencil $(W, W') \subset |D|$, containing $W$, such that
  $W'_W = E_0 + E_1$, with $E_0 \geq 0$ and $E_1 \in \Lambda$. By Hironaka's
  theory 
\cite{Hir}, we may use blowing-ups with smooth centers contained in $W \cap W'$, such that after taking the proper transforms of the elements of our pencil, to get:\begin{enumerate}
\item[(a)] $  {\rm Supp} (E_0)$ has normal crossing;
\item[(b)] $\Lambda$ is base points free (so $\varphi : W \to \pn {n-1}$ is a birational morphism).
\end{enumerate}
Further blowing-up of the components of ${\rm Supp} (E_0)$ allows to assume $E_0=0$ so $D_W$ is linearly equivalent to $E$. Using the previous exact
sequence and the fact that $H^1(X, \sO _X) = 0$, it follows that ${\rm Bs}|D| =\emptyset$. Finally, $D^n = (D_W)_W^{n-1} = 1$, so $\varphi$ is a
birational morphism to $\pn n$.
\qed

\begin{example*}\label{sad} (\cite[Example 1.4]{IoNa}) Let $X \subset \pn {n+d-2}$ be a non-degenerate projective
    variety of dimension $n\geq 2$ and degree $d\geq 3$, which   is not a
    cone.
Then $X$ is rational, unless it is a smooth cubic hypersurface, $n\geq 3$.
If $X$ is singular,  by projecting  from a singular  point we get a variety of
minimal degree, birational to $X$. So $X$ is rational. If $X$ is not linearly
normal, $X$ is isomorphic to a variety of minimal degree. Hence we may assume
$X$ to be smooth and linearly normal. One sees easily that such a linearly
normal, non-degenerate manifold $X\subset \pn {n+d-2}$ has anticanonical
divisor linearly equivalent to $n-1$ times the hyperplane section, i.e., they
are exactly the so-called  ``classical del Pezzo manifolds''. They were
classified  by Fujita in a series of papers; see \cite{FuBook} or \cite{IP}. 
Independently of their classification, the following  simple argument shows that such manifolds are rational when $d \ge 4$.
 Consider the surface $W$ obtained by intersecting $X$ with $n-2$
general hyperplanes. Note that $W$ is a non-degenerate, linearly normal surface 
of degree $d$ in $\pn d$, so it is a del Pezzo surface. As such, $W$ is known
to admit a representation $\varphi:W \to \pn 2$ as the blowing-up of $9-d$
points. Let $L\subset W$ be the pull-back via $\varphi$
of a general line in $\pn 2$. It is easy to see that $L$ is a cubic rational
curve in the embedding of $W$ into $\pn d$. So, for $d\geq 4$, $L$ is contained
in a hyperplane of $\pn d$. This shows that the conditions of the
Theorem~\ref{trei-sase} are fulfilled for $X$, $|D|$ being the system of
hyperplane sections. We also see that Theorem~\ref{trei-sase} is sharp, as the
previous argument fails exactly for the case of cubics. 
\end{example*}

\begin{rem*}\label{closing}
In closing,  we mention three possible   generalizations of the problem of extending morphisms from  ample divisors on $X$. 

\begin{enumerate}

\item[(1)] The smoothness assumption on $X$ may be relaxed by allowing normal singularities.  Let  $Y$  be a  smooth divisor in $X$ ($X$ is smooth), and let us  only suppose  that $Y$ has ample  normal bundle. Then  a well-known result (\cite{Hrs}) shows that there is a birational map 
$\psi:X\to X'$, which is an isomorphism along $Y$, such that $\psi(Y):=Y'\subset X'$ is ample and $X'$ is normal. 
See e.g.,  \cite{BadHD,BadSing} and \cite{CaFl} for results in this direction.

\item[(2)]   Consider a smooth section $Y\subset X$ of the appropriate expected dimension $n-{\rm rk}\sE$ of an ample vector bundle $\sE$ on an $n$-fold  $X$. Note that  a Lefschetz type theorem for ample vector bundles,  due to Sommese \cite{SoLef}, implies that  the restriction  to $Y$ gives  an isomorphism 
${\rm Pic}(X)\cong {\rm Pic}(Y)$.
See e.g., \cite{AOInt,dFL,LaMa,AO2}  and \cite{anoc} for results of this type.

\item[(3)]
In the same spirit, let us consider a smooth subvariety $Y$ of a manifold $X$ such that ${\rm codim}_XY\geq 2$, and $Y$ has ample normal bundle.  Further, let us  add the Lefschetz type assumption that ${\rm Pic}(Y)\cong {\rm Pic}(X)$. 
Then  one can study extensions of rationally connected fibrations $p:Y\to Z$ onto a normal projective  variety $Z$. See
 \cite{BdFL} and \cite{Occ}  for results in this direction.

\item[(4)]  The very recent paper \cite{Wa} classifies pairs $(X,Y)$, when $Y\subset X$ is an ample divisor which is a homogeneous manifold.
\end{enumerate}

\end{rem*}

\section*{Acknowledgments}
 
 We thank the referee for several useful comments.

\bibliographystyle{amsalpha}





\end{document}